\documentclass[11pt]{article}
\usepackage{lmodern}
\usepackage[T1]{fontenc}
\usepackage{hyperref}
\usepackage{amsfonts,amsmath,amssymb,amsopn}
\usepackage[all]{xy}
\usepackage{vmargin}
\usepackage{makeidx}
\title{Constructible $\nabla$-modules on curves}
\author{Bernard Le Stum\thanks{bernard.le-stum@univ-rennes1.fr}}

\date{Version of \today}
\newtheorem{thm}{Theorem}[section]
\newtheorem{prop}[thm] {Proposition}
\newtheorem{cor}[thm] {Corollary}
\newtheorem{lem}[thm] {Lemma}
\newtheorem{dfn}[thm] {Definition}
\begin{document}

\maketitle

\begin{center}
\textbf{Abstract}
\end{center}

Let $\mathcal V$ be a discrete valuation ring of mixed characteristic with perfect residue field.
Let $X$ be a geometrically connected smooth proper curve over $\mathcal V$.
We introduce the notion of constructible convergent $\nabla$-module on the analytification  $X_{K}^{\mathrm{an}}$ of the generic fibre of $X$.
A constructible module is an $\mathcal O_{X_{K}^{\mathrm{an}}}$-module which is not necessarily coherent, but becomes coherent on a stratification by locally closed subsets of the special fiber $X_{k}$ of $X$.
 The notions of connection, of (over-) convergence and of Frobenius structure carry over to this situation.
We describe a specialization functor from the category of constructible convergent $\nabla$-modules to the category of $\mathcal D^\dagger_{\hat X \mathbf Q}$-modules.
We show that if $X$ is endowed with a lifting of the absolute Frobenius of $X$, then specialization induces an equivalence between constructible $F$-$\nabla$-modules and perverse holonomic $F$-$\mathcal D^\dagger_{\hat X \mathbf Q}$-modules.

\tableofcontents

\section*{Introduction}

When we look for a category of coefficients for a cohomological theory, we usually get one which is too small inside another one which is too big and we aim at the perfect one that will sit in between.

For example, if we are interested in the singular cohomology of an algebraic variety over $\mathbf C$, we have on one side the category of local systems of finite dimensional vector spaces and on the other side the category of all sheaves of vector spaces.
The perfect category will be the category of constructible sheaves.
Similarly, for de Rham cohomology, we have on one side the category of coherent modules endowed with an integrable connection and on the other side the category of all $\mathcal D$-modules.
The perfect one will be the category of regular holonomic $\mathcal D$-modules.
Last, we may consider the category of finitely presented crystals and the category of all modules on the infinitesimal site.
In between, we have the category of constructible pro-crystals (unpublished note of Deligne (\cite{Deligne*}). 
We have equivalences of categories between local systems, modules with an integrable connection and finitely presented crystals (we need to add some regularity conditions when moving from analytic to algebraic side).
At the derived category level, we also have an equivalences between constructible sheaves, regular holonomic $\mathcal D$-modules and constructible pro-crystals (\cite{Kashiwara84} and \cite{Deligne*}).

If we are interested in the $p$-adic cohomology of a variety of characteristic $p > 0$, there is no analog to the first theory.
The closest would be the \'etale $p$-adic cohomology which is not satisfying.
There is a very good analog to the second theory, which is the theory of arithmetic $\mathcal D$-modules (and overconvergent isocrystals) developed by Berthelot, Virrion, Huyghe, Trihan, Caro and others (see \cite{Berthelot02} for an overview).
More recently, I started to develop a crystalline theory in \cite{LeStum11}.
I showed that the category of overconvergent isocrystals is equivalent to the category of finitely presented crystals on the overconvergent site.
When I was visiting the university of Tokyo in october 2009, Shiho suggested that we should look for the perfect category in this new theory, and prove that it is equivalent to the one developed by Berthelot (an \emph{overconvergent Deligne-Kashiwara correspondence}).
I had a very naive idea of what this category of constructible overconvergent crystals should be and gave a lecture on this topic in March 2010 at Oxford University.
The aim of this article is to show that the strategy works for curves.

More precisely, if $X$ is a smooth proper curve over a complete discrete valuation ring $\mathcal V$ of mixed characteristic $p$ with fraction field $K$ and perfect residue field $k$, we introduce the notion of \emph{constructible convergent $\nabla$-module} $E$ on $X_{K}^{\mathrm{an}}$ and show that its specialization $\mathrm R\widetilde{\mathrm{sp}}_{*}E$ is a \emph{perverse} (complex of) $\mathcal D^\dagger_{\hat X\mathbf Q}$-module.
When there is a Frobenius on $X$, we show that specialization induces an equivalence between \emph{constructible $F$-$\nabla$-modules} on $X_{K}^{\mathrm{an}}$ and perverse holonomic $F$-$\mathcal D^\dagger_{\hat X\mathbf Q}$-modules.
It should be remarked that the overconvergent isocrystals on an open subset $U$ of $X_{k}$ form a full subcategory of the category of all constructible convergent $\nabla$-modules on $X_{K}^{\mathrm{an}}$ and that our theorems extend Berthelot's equivalence theorems for specialization of (over-) convergent isocrystals.

Our definitions are very general and apply to any formal $\mathcal V$-scheme.
A constructible overconvergent $\nabla$-module is an (not necessarily coherent) $\mathcal O_{X_{K}^{\mathrm{an}}}$-module $E$ with a convergent (automatically integrable here) connection and we require that there exists a finite covering of the special fibre $X_{k}$ by locally closed subset $Y$ such that, if $i_{Y} : ]Y[ \hookrightarrow X_{K}^{\mathrm{an}}$ denotes the inclusion map, then $i_{Y}^{-1}E$ is a coherent $i_{Y}^{-1}\mathcal O_{X_{K}^{\mathrm{an}}}$-module.
In this definition, we view $X_{K}^{\mathrm{an}}$ as a Berkovich analytic space but we define $\mathrm R\widetilde{\mathrm{sp}}_{*}E := \mathrm R \mathrm{sp}_{*}E_{G}$ as the derived specialization for the Tate topology.

We can give an explicit description of constructible convergent $\nabla$-modules and their specialization.
If $D$ is a smooth relative divisor on $X$ with affine open complement $\mathrm{Spec}\; A$, we can consider $A^\dagger_{K}$ which is the generic fiber of the weak completion of $A$, and for each $a \in D_{K}$, the Robba ring $\mathcal R_{a}$ on $K(a)$.
Then, a constructible convergent $\nabla$-module is equivalent to the following data for some sufficiently big divisor $D$:  an overconvergent $\nabla$-module $M$ (of finite type) over $A^\dagger_{K}$ and for each $a \in D_{K}$, a finite dimensional $K(a)$-vector space $H_{a}$ and a horizontal $A^\dagger_{K}$-linear map
$$
M \to \mathcal R_{a} \otimes_{K(a)} H_{a}.
$$
Now, we can describe $\mathrm R\widetilde{\mathrm{sp}}_{*}E$ as follows.
Let $\mathcal U' = \mathrm{Spf} \mathcal A'$ be an affine open subset of $\hat X$.
First of all, if $\mathcal D' = \mathcal U' \cap \hat D$ is defined in $\mathcal U'$ by an equation $f$, we may consider the weak completion $\mathcal A'[1/f]^\dagger$.
Also, we can always write each $\mathcal R_{a}$ as an extension of a $K(a)$-vector space $\delta_{a}$ by the ring $\mathcal O_{a}^{\mathrm{an}}$ of convergent functions on the open unit disc.
We extend scalars on the left, project on the right and add all the above maps in order to get the complex
$$
\mathcal A'[1/f]^\dagger \otimes_{A^\dagger_{K}} M \to \bigoplus_{a \in \mathcal D'_{K}} \delta_{a} \otimes_{K(a)} H_{a}.
$$
This is $\mathrm R\widetilde{\mathrm{sp}}_{*}E$ on $\mathcal U'$.
It is naturally a complex of $\mathcal D^\dagger_{\hat X\mathbf Q}$-modules on $\mathcal U'$.
It is perverse in the sense that it has flat cohomology in degree 0, finite support in degree 1 and no other cohomology.

As already mentioned, we use Berkovich theory of ultrametric geometry (see \cite{Berkovich90} and \cite{Berkovich93} or \cite{Ducros07}) and will try to recall the basic constructions when we meet them.
We will generally avoid reference to basic results proved in the context of rigid analytic geometry and rather reprove them here in this new language when necessary.
Anyway, the reader should keep in mind that a significant part of what follows is merely a reorganization of some of Berthelot's material using Berkovich theory instead of Tate's.

Note that I will use the terminology \emph{finite module} for finitely presented module and \emph{finite torsion module} if we ever have to consider a module which is finite as set.
I will call {$\nabla$-module} a module endowed with a (integrable) connection.
We denote with an upper $*$ the inverse image for a morphism of ringed spaces when the upper $-1$ is used for inverse image for a morphism of topological spaces.
Also, we will only consider real numbers $\eta$ inside $\sqrt{|K|}$ where $K$ is our fixed ultrametric field.
Finally, I will systematically identify a sheaf on a topological space reduced to one point with its global sections.

Many thanks to Daniel Caro and Pierre Berthelot who helped me a lot when I was struggling with arithmetic $\mathcal D$-modules problems.
I am indebted too to Florian Ivorra who helped me clarify some questions related to adjunction of derived functors.
I also recall that the problem which is solved here was pointed out by Atsushi Shiho.

\section{Constructible modules}

We fix a complete ultrametric field $K$ with ring of integers $\mathcal V$, maximal ideal $\mathfrak m$ and residue field $k$.
We will assume that $K$ is non trivial and fix some $\pi \in K$ with $0 < |\pi| < 1$.
We will also assume that $K$ has characteristic zero, that $k$ has positive characteristic $p$, that the valuation is discrete and that $k$ is perfect.
When Frobenius enters the game, we also fix an isometry $\sigma$ on $K$ that lifts the absolute Frobenius of $k$.

We let $X$ be a geometrically connected proper smooth $\mathcal V$-scheme of relative dimension one.

We may consider its special fiber $X_{k}$ as well as its generic fiber $X_{K}$ which are proper smooth geometrically connected curves and we will be particularly interested in the \emph{Berkovich} analytification $X_{K}^{\mathrm{an}}$ of $X_{K}$.
For the moment, we will simply use the fact that, since $X$ is proper, we can identify $X_{K}^{\mathrm{an}}$ with the generic fiber $\hat X_{K}$ of the completion of $X$.
In particular, $X_{K}^{\mathrm{an}}$ is covered (for Tate topology) by the affinoid domains $\mathcal U_{K} = \mathcal M(\mathcal A_{K})$ (set of continuous multiplicative semi-norms) if $\hat X$ is covered by some affine subsets $\mathcal U = \mathrm{Spf}\; \mathcal A$ (set of open primes).

We will consider the specialization map
$$
\mathrm{sp} : X_{K}^{\mathrm{an}} = \widehat X_{K} \to \widehat X \simeq X_{k}
$$
(where the last morphism is a homeomorphism).
Locally, the specialization of a point $x \in \mathcal U_{K}$ is the open prime
$$
\mathfrak p = \{f \in \mathcal A, \quad |f(x)| < 1\} \in \mathrm{Spf}\; \mathcal A = \mathcal U
$$
(see section 1 of  \cite{Berkovich94} for the general construction).
Note that if $x$ is a semi-norm on a ring and $f$ a function in this ring, we write $|f(x)|$ instead of $x(f)$.

When $Y$ is a subset of $X_{k}$, we will consider its \emph{tube}
$$
]Y[ := \mathrm{sp}^{-1}(Y) \subset X_{K}^{\mathrm{an}}.
$$
This is mostly used when $Y$ is locally closed, in which case $]Y[$ is an analytic domain in $X_{K}^{\mathrm{an}}$.
We will give an explicit description below, but it should first be remarked that, since $X_{k}$ is a connected curve, any locally closed subset $Y$ of $X_{k}$ is necessarily open or closed.
Moreover, any open subset $U \subset \neq X_{k}$ is affine.
Also, any closed subset $Z \subset \neq X_{k}$ is a finite set of closed points.

When $U$ is open in $X_{k}$, then $]U[$ is a \emph{closed} subset of $X_{K}^{\mathrm{an}}$.
More precisely, if $\mathcal U$ is the formal lifting of $U$, then  $]U[ = \mathcal U_{K}$ (which is affinoid and therefore compact).
Conversely, if $Z$ is a closed subset of $X_{k}$, then $]Z[$ is an \emph{open} subset of $X_{K}^{\mathrm{an}}$: more precisely, if $Z$ is defined inside some formal affine open subset $\mathcal U$ by an equation $f = 0 \mod \mathfrak m$, then $]Z[$ is defined in $\mathcal U_{K}$ by $|f(x)| < 1$ (being open is local for Tate topology).
It will also be convenient to consider the tubes $[Z]_{\eta}$ and $]Z[_{\eta}$ of radius $\eta < 1$ defined by $|f(x)| \leq \eta$ and $|f(x)| < \eta$ respectively (they are independent of the choice of $f$ when $\eta$ is close to $1$ ).

The specialization map $\mathrm{sp} : X_{K}^{\mathrm{an}} \to X_{k}$ is surjective and may be used to give a description of $X_{K}^{\mathrm{an}}$.
First of all, there is only one point above the generic point $\xi$ of $X_k$ that we shall therefore denote by $]\xi[$ and call the \emph{generic point} of $X_{K}^{\mathrm{an}}$ even if it depends on the model.
Actually, if $U \subset X_{k}$ is a non empty affine subset, then $]\xi[$ is the \emph{Shilov boundary} of $]U[$ (the point where any function reaches its maximum).
It is also the usual boundary of $]U[$ as subset of the topological space $X_{K}^{\mathrm{an}}$ and the absolute boundary of $]U[$ in Berkovich sense.
For example, in the case of $X = \mathbf P^1_{\mathcal V}$, then $]\xi[$ is the Gauss norm.

By definition, $X_K^{\mathrm{an}} \setminus ]\xi[$ is the disjoint union of all open subsets $]x[$ for $x$ a closed point in $X_k$.
These open subsets $]x[$ are called the \emph{residual classes}.
In other words, the generic point is used to connect the residual classes together (recall that $X_{K}^{\mathrm{an}}$ is simply connected).
We can also mention that $\cap ]U[ = ]\xi[$ when $U$ runs trough all the non empty (affine) open subsets of $X_{k}$.

Before giving an explicit description of the residual classes, recall that there exists also a canonical map $X_{K}^{\mathrm{an}} \to X_{K}$.
If $\mathrm{Spec}\; A \subset X$ is an open subset, the analytification of $\mathrm{Spec}\; A_{K}$ is the set $\mathcal M^{\mathrm{alg}}(A_{K})$ of all multiplicative semi-norms $x$ on $A$ and $x$ is sent to its kernel
$$
\mathfrak p = \{f \in A, \quad |f(x)| = 0\} \in \mathrm{Spec}\; A_{K}.
$$
The map $X_{K}^{\mathrm{an}} \to X_{K}$ induces a bijection between the rigid points of $X_{K}^{\mathrm{an}}$ (whose residue fields are finite over $K$) and the closed points of $X_{K}$
(note that all the non rigid points are then sent to the generic point of $X_{K}$).
We will implicitly identify the closed points of $X_{K}$ with the rigid point of $X_{K}^{\mathrm{an}}$.

We call a closed point $a \in X_{K}$ a lifting of a closed point $x \in X_{k}$ if $\mathrm{sp}(a) = x$.
Note that $K(a)$ is a finite extension of $K$ and inherits a structure of complete ultrametric field whose ring of integers will be denoted $\mathcal V(a)$.
The point $a$ is said to be \emph{unramified} if $K(a)/K$ is an unramified extension (i.e $\mathcal V(a)$ is \'etale over $\mathcal V$).

Here is a standard result on residual classes:

\begin{prop} \label{pointdisc}
If $x \in X_{k}$ is a closed point and $a$ is an unramified lifting of $x$, then $]x[ \simeq \mathbf D_{K(a)}(0, 1^-)$.
Actually, we have $[x]_{\eta} \simeq \mathbf D_{K(a)}(0, \eta^+)$ for $\eta < 1$ and close to $1$.
\end{prop}

\textbf {Proof: }
This follows from Berthelot's weak fibration theorem (see lemma 4.4 of \cite{Berkovich99}) as we can now recall.
When $x$ is a rational point, we can write $x$ as the unique zero of an \'etale map
$$
t : \mathcal U = \mathrm{Spf}\; \mathcal A \to \widehat {\mathbf A}^1_{\mathcal V}.
$$
The corresponding morphism $\mathcal U_{K} \to \mathbf D_{K}(0, 1^+)$ will induce an isomorphism $[x]_{\eta} \simeq \mathbf D(0, \eta^+)$.
In general, we first have to extend the basis and consider the projection $X_{\mathcal V(a)} \to X$.
We pick up some point $y$ over $x$.
Then, the corresponding morphism $X_{K(a)}^{\mathrm{an}} \to X_{K}^{\mathrm{an}}$ will induce an isomorphism $[y]_{\eta} \simeq [x]_{\eta}$.
$\quad \Box$

We are finished with the study of finite closed subsets and we consider now the case of a non empty affine open subset $U \subset X_{k}$.
We may always lift  \emph{non canonically} $U$ to an affine open subset $\mathrm{Spec}\; A \subset X$.
Since $k$ is perfect, we may even assume that $\mathrm{Spec}\; A$ is the open complement of (the support of) a smooth (relative) divisor $D \subset X$.
In other words, the complement of $U$ in $X_{k}$ is $D_{k} =: \{x_{1}, \ldots, x_{n}\}$ and the complement of $\mathrm{Spec}\; A_{K}$ in $X_{K}$ is $D_{K} = \{a_{1}, \ldots, a_{n}\}$ where each $a_{i}$ is an unramified lifting of $x_{i}$.
In this case, choosing $A$ amounts to fixing the ``center'' $a_{i}$ of each disc $]x_{i}[$.

Actually, we will only use the weak completion $A^\dagger$ of $A$ and more precisely its generic fiber $A^\dagger_K$.
In fact, if $Z$ denotes the closed complement of $U$, the subsets $V_{\lambda} := X_{K}^{\mathrm{an}} \setminus ]Z[_{\lambda}$ form a cofinal family of affinoid neighborhoods of $]U[$ and $A^\dagger_{K} = \varinjlim A_{\lambda}$ if we write $V_{\lambda} =: \mathcal M(A_{\lambda})$.
In particular, the algebra $A^\dagger_{K}$ does not depend on the algebraic lifting but only on $U$.
However, the geometric construction can be useful and we will stick to the original definition. 

It is also easy to see that the subsets $V_{\lambda}$ for various $\lambda$ \emph{and} $U$ form a cofinal family of affinoid neighborhoods of the generic point $]\xi[$.
Actually, since $X_{K}^{\mathrm{an}}$ is compact, any open neighborhood $V$ of $]\xi[$ has a compact complement $C \subset \cup ]x[ = \cup ]x[_{\lambda}$ and therefore, $C \subset \cup_{\mathrm{finite}} ]x[_{\lambda} = ]T[_{\lambda}$.
It follows that $V_{\lambda} := X_{K}^{\mathrm{an}} \setminus ]T[_{\lambda} \subset V$.
In particular, we have
$$
\mathcal O_{X_{K}^{\mathrm{an}}, ]\xi[} = \varinjlim A^\dagger_{K}
$$
when $\mathrm{Spec}\; A$ runs through the non empty affine open subsets of $X$.
Note that $\mathcal O_{X_{K}^{\mathrm{an}}, ]\xi[}$ is a henselian field and that $|f(]\xi[)| = \|f\|_{]U[}$ if $f \in A^\dagger_{K}$ and $\mathrm{Spec}\; A$ is a lifting of $U$.

If $Y$ is a subset of $X_{k}$, we will denote by
$$
i_{Y} : ]Y[  \hookrightarrow X_{K}^{\mathrm{an}}
$$
the corresponding inclusion.
Recall that $i_{U}$ is a \emph{closed} embedding when $U$ is an open subset in $X_{k}$ and that $i_{Z}$ is an \emph{open} embedding when $Z$ is a closed subset of $X_{k}$.

\begin{dfn}
An $\mathcal O_{X_{K}^{\mathrm{an}}}$-module $E$ is \emph{constructible} (resp. {constructible-free}) if there exists a finite covering of $X_{k}$ by locally closed subsets $Y$ such that $i_{Y}^{-1} E$ is a coherent (resp. free) $i_{Y}^{-1} \mathcal O_{X_{K}^{\mathrm{an}}}$-module.
\end{dfn}

This definition makes sense in a very general situation but we really want to stick to curves here in order to give an explicit description of these constructible modules.
Nevertheless, the following is purely formal:

\begin{prop}
Constructible (resp. constructible-free) modules form an abelian (resp. additive) subcategory stable under internal Hom, extensions and tensor product.
\end{prop}

\textbf {Proof: }
Since the restriction maps $E \mapsto i_Y^{-1}E$ are exact and commute with internal hom and tensor product, all these assertions follow from the analogous results on coherent (resp. free) modules.
$\quad \Box$

Note that if $Y' \subset Y$ and  $i_{Y}^{-1} E$ is coherent, free or locally free, so is $i_{Y'}^{-1} E$ by restriction.
It should also be remarked that if $Z \subset X_{k}$ is a closed subset, then $i_{Z}^{-1} E = E_{|]Z[}$ is the usual restriction to an open subset.
However, when $U \subset X_{k}$ is open, then $i_{U}$ is a closed embedding and $i_{U}^{-1} E$ is the sheaf of sections defined in the neighborhood of the compact subset $]U[$.
More on this soon.

At some point, we will need to use some theorems that were formulated in the language of rigid analytic geometry.
It will therefore be necessary to have a dictionary at our disposal.
We want to explain this right now.
We may consider the Grothendieck-Tate topology on $X_{K}^{\mathrm{an}}$: admissible open subsets are analytic domains and coverings are \emph{Tate} covering (see for example \cite{Berkovich93}, page 25-26).
If $V$ is an analytic domain of $X_{K}^{\mathrm{an}}$, we will denote by $V_{G}$ the associated Grothendieck space.
There is an obvious continuous map of Grothendieck spaces $\pi_{V} : V_{G} \to V$ which is simply the identity on underlying sets.
Actually, $\pi_{V}^*$ is exact and fully faithful and induces an equivalence on coherent sheaves.
We will write $E_{G} := \pi_{V}^{*}E$.
Finally, the set  $V^{\mathrm{rig}}$ of rigid points of $V$ has a natural structure of rigid analytic variety and the topos of $V^{\mathrm{rig}}$ is equivalent to the topos of $V_{G}$ (affinoid domains and affinoid coverings coincide).
We will denote by $E_{0}$ the rigid sheaf corresponding to $E_{G}$.

The next comparison theorem was proved in corollary 2.2.13 of \cite{LeStum11} but it will not be long to recall how it works.

\begin{prop} \label{comprig}
If $U$ is an open subset of $X_{k}$, the functor $E \mapsto (i_{U*}E)_{0}$ induces an equivalence between the category of coherent $i_{U}^{-1}\mathcal O_{X_{K}^{\mathrm{an}}}$-modules on $]U[$ and the category of coherent $j^\dagger \mathcal O_{X_{K}^{\mathrm{rig}}}$-modules.
\end{prop}

\textbf {Proof: }
If $Z$ is the closed complement of $U$ and $\lambda < 1$, we consider the affinoid domain $V_{\lambda} = X_{K}^{\mathrm{an}} \setminus ]Z[_{\lambda}$.
Since $X_{K}^{\mathrm{an}}$ is locally compact, the category of coherent $i_{U}^{-1}\mathcal O_{X_{K}^{\mathrm{an}}}$-modules is equivalent to the direct limit of the category of coherent $\mathcal O_{V_{\lambda}}$-modules (see corollary 2.2.5 of \cite{LeStum11} for example).
And we also know that the category of coherent $j^\dagger \mathcal O_{X_{K}^{\mathrm{rig}}}$-modules is equivalent to the direct limit of the category of coherent $\mathcal O_{V^{\mathrm{rig}}_{\lambda}}$-modules (theorem 5.4.4 of \cite{LeStum07}).
$\quad \Box$

In general, we will not use this comparison theorem because most results are a lot simpler to prove in the Berkovich topology and it sounds unnatural to use the rigid topology in order to derive them.
For example, we have the following (so called theorem A and B) which is essentially proposition 5.4.8 of \cite{LeStum07}:

\begin{prop} \label{dageq}
If $U$ is an affine open subset of $X_{k}$, then $i_{U}^{-1} \mathcal O_{X_{K}^{\mathrm{an}}}$ is a coherent ring.
Moreover, if $\mathrm{Spec}\; A \subset X$ is an algebraic lifting of $U$, then
$$
\mathrm R \Gamma(]U[, i_{U}^{-1} \mathcal O_{X_{K}^{\mathrm{an}}} ) = A^\dagger_{K}.
$$
Finally, the functor $\mathrm R \Gamma(]U[, -)$ induces an equivalence between coherent (resp. locally free, resp. free) $i_{U}^{-1} \mathcal O_{X_{K}^{\mathrm{an}}}$-modules and finite (resp. finite projective, resp. finite free) $A^\dagger_{K}$-modules.
\end{prop}

\textbf {Proof: }
Since $]U[$ is a compact subset of $X_{K}^{\mathrm{an}}$ which is locally compact and $ \mathcal O_{X_{K}^{\mathrm{an}}}$ is coherent, the first assertion is purely formal (see for example proposition 2.2.3 of \cite{LeStum11}).
For the same reason, any coherent $i_{U}^{-1} \mathcal O_{X_{K}^{\mathrm{an}}}$-module $E$ is the restriction of a coherent module on some neighborhood $V$ of $]U[$ in $X_{K}^{\mathrm{an}}$ and we get an equivalence of categories if we allow the shrinking of $V$ (see corollary 2.2.5 of \cite{LeStum11}).
As mentioned above, there exists a cofinal family of affinoid neighborhoods $V_{\lambda} = \mathcal M (A_{\lambda})$ with $A^\dagger_{K} = \varinjlim A_{\lambda}$.
The rest of the proposition follows then immediately from the analogous results on affinoid varieties and the properties of filtered direct limits.
$\quad \Box$

\begin{lem}\label{trivlem}
If $E$ is a constructible (resp. constructible free) module and $Z \subset X_{k}$ is a finite closed subset, then $i_{Z}^{-1}E$ is coherent (resp. locally free).
\end{lem}

\textbf {Proof: }
We may assume that $Z$ is reduced to one point $x$.
Then, there exists a locally closed subset $Y \subset X_{k}$ with $x \in Y$ and $i_{Y}^{-1} E$ coherent (resp. free) on $]Y[$.
We may then restrict to $]x[$.
$\quad \Box$

\begin{prop}
If $E$ is a constructible module, there exists a non-empty open subset $U \subset X$ such that $i_{U}^{-1}E$ is a free $i_{U}^{-1} \mathcal O_{X_{K}^{\mathrm{an}}}$-module.
\end{prop}

\textbf {Proof: }
By definition, if $\xi$ is the generic point of $X_k$, there exists a locally closed subset $U \subset X_k$, which is necessarily open, such that $\xi \in U$ and $i_{U}^{-1} E$ is a coherent $i_{U}^{-1} \mathcal O_{X_{K}^{\mathrm{an}}}$-module.
It follows that the stalk $E_{]\xi[}$ of $E$ at $]\xi[$ is a finite dimensional vector space (recall that $\mathcal O_{X_{K}^{\mathrm{an}},]\xi[}$ is a field).
At this point, it is convenient to write $i_{U}^{-1} E = i^{-1}\mathcal F$ where $i : ]U[ \hookrightarrow V$ if the inclusion of $]U[$ in some neighborhood and $\mathcal F$ is a coherent $\mathcal O_{V}$-module.
Since $\mathcal F_{]\xi[} = E_{]\xi[}$ is a finite dimensional vector space and $\mathcal F$ is coherent, we know that there exists a neighborhood of $\xi$ in $V$ on which $\mathcal F$ becomes free.
We saw above that the subsets $V_{\lambda} = X_{K}^{\mathrm{an}} \setminus ]T[_{\lambda}$ with $T$ finite closed and $\lambda < 1$ form a cofinal family of neighborhoods of $]\xi[$.
Therefore, there exists such a $T$ and $\lambda$ with $\mathcal F_{|V_{\lambda}}$ free.
We may then remove a finite number of closed points in $U$ and assume that $U \cap T = \emptyset$ in order to get $i_{U}^{-1} E$ free.
$\quad \Box$

\begin{cor} \label{freeness}
Let $E$ be a constructible module.
Then, the following are equivalent:
\begin{enumerate}
\item $E$ is constructible-free
\item $i_{Z}^{-1}E$ is locally free whenever $Z \subset X_{k}$ is a finite closed subset.
\item $i_{x}^{-1}E$ is free whenever $x \in X_{k}$ is a closed point
\end{enumerate}
\end{cor}

\textbf {Proof: }
We know from proposition \ref{pointdisc} that if $x \in X_{k}$ is any closed point, then $]x[$ is a disc.
Moreover, since the valuation is discrete, a coherent locally free module on a disk is automatically free.
It follows that assertion 2 and 3 are equivalent.
And we showed in lemma \ref{trivlem} that they are automatically fulfilled when $E$ is constructible-free.
Conversely, since $E$ is constructible, we know from the proposition that there exists a non-empty open subset $U \subset X_{k}$ such that $i_{U}^{-1}E$ is free.
If we assume that assertion 3 holds, then, in particular, for all $x \not \in U$, $i_{x}^{-1}E$ is free.
Now, $U$ and the points $x \not \in U$ form a finite covering of $X_{k}$ by locally closed subsets.
$\quad \Box$

\begin{cor}
An $\mathcal O_{X_{K}^{\mathrm{an}}}$-module $E$ is constructible-free if and only if there exists a finite covering of $X_{k}$ by locally closed subsets $Y$ such that $i_{Y}^{-1} E$ is a locally free $i_{Y}^{-1} \mathcal O_{X_{K}^{\mathrm{an}}}$-module.
\end{cor}

\textbf {Proof: }
The condition is clearly necessary.
Conversely, any such module is constructible and we may apply the previous corollary.
$\quad \Box$

\begin{cor} \label{locdec}
An $\mathcal O_{X_{K}^{\mathrm{an}}}$-module $E$ is constructible (resp. constructible-free) if and only if there exists a non-empty affine open subset $U \subset X_{k}$ with closed complement $Z$ such that $i_{U}^{-1} E$ is a coherent (resp. locally free) $i_{U}^{-1} \mathcal O_{X_{K}^{\mathrm{an}}}$-module and $i_{Z}^{-1}E$ is coherent (resp. locally free).
Moreover, we may assume that $i_{U}^{-1} E$ is free.
\end{cor}

\textbf {Proof: }
It follows from the previous corollary that the condition is necessary.
Conversely, if $E$ is constructible, we can find such a $U$ thanks to the proposition and then apply lemma \ref{trivlem} to the complement $Z$ of $U$ in $X_{k}$.
$\quad \Box$

\begin{cor} \label{coht}
Let $E$ be a constructible-free module and $Y \subset X_{k}$ a locally closed subset.
If $i_{Y}^{-1}E$ is coherent, it is necessarily locally free.
\end{cor}

\textbf {Proof: }
If $a \in ]Y[$ specializes to $x \in X_{k}$, then $a \in ]x[$ and we know that $i_{x}^{-1}E$ is free.
It follows that the stalk $E_{a}$ of $E$ at $a$ is free.
A coherent module whose stalks are free is necessarily locally free.
$\quad \Box$

If $Z \subset X_k$ is a closed subset, then $i_Z : ]Z[ \hookrightarrow X_{K}^{\mathrm{an}}$ is an open immersion and we may consider the extension by zero $i_{Z!}$ outside $Z$ which is a left adjoint functor to $i_Z^{-1}$.
 
\begin{prop} \label{excons}
An $\mathcal O_{X_{K}^{\mathrm{an}}}$-module $E$ is constructible (resp. constructible free) if and only if there exists an exact sequence
$$
0 \to i_{Z!}E_{Z} \to E \to i_{U*} E_{U} \to 0 
$$
where $U$ is a non-empty affine open subset of $X_{k}$ with closed complement $Z$, $E_{U}$ is a coherent (resp. locally free) $i_{U}^{-1} \mathcal O_{X_{K}^{\mathrm{an}}}$-module and $E_{Z}$ is coherent (resp. locally free).
Moreover, we may assume that $i_{U}^{-1} E$ is free.
\end{prop}

\textbf {Proof: }
If we are given such and exact sequence, we may pull back along $i_{U}$ and $i_{Z}$ in order to obtain $E_{U} \simeq i_{U}^{-1}E$ and $E_{Z} \simeq i_{Z}^{-1}E$.
And conversely, using corollary \ref{locdec}, we can find a non-empty affine open subset $U \subset X_{k}$ with closed complement $Z$ such that $i_{U}^{-1} E$ is a coherent (resp. locally free) $i_{U}^{-1} \mathcal O_{X_{K}^{\mathrm{an}}}$-module and $i_{Z}^{-1}E$ is coherent (resp. locally free).
There is an exact sequence
$$
0 \to i_{Z!}i_{Z}^{-1}E \to E \to i_{U*} i_{U}^{-1}E \to 0 
$$
and we can set $E_{U} :=  i_{U}^{-1}E$ and $E_{Z} := i_{Z}^{-1}E$.
$\quad \Box$

The notion of constructibility is stable under pull back by the very definition but we will only use the following particular case:

\begin{prop}
Assume that $F: X \to X$ is a $\sigma$-linear lifting of the absolute Frobenius of $X_{k}$.
If $E$ is a constructible module on $X_{K}^{\mathrm{an}}$, then $F_{K}^{\mathrm{an}*}E$ is also constructible.
\end{prop}

\textbf {Proof: }
Follows from the fact that coherence is stable under pull-back. 
$\quad \Box$

\section{Classification}

Recall that if $Y$ is a locally closed subset of $X_{k}$, we denote by $i_{Y} : ]Y[ \hookrightarrow X_{K}^{\mathrm{an}}$ the inclusion map.
Recall also that $\xi$ denotes the generic point of $X_{k}$.

\begin{dfn} \label{robdir}
Let $T \subset X_{k}$ be a non empty finite closed subset.
If $\mathcal F$ is a coherent $\mathcal O_{]T[}$-module, we will write $\mathcal O^{\mathrm{an}}(\mathcal F) := \Gamma(]T[, \mathcal F)$.
The \emph{Robba module} of $\mathcal F$ is
$$
\mathcal R(\mathcal F) := (i_{T*}\mathcal F)_{]\xi[}.
$$
and the \emph{Dirac space} of $\mathcal F$ is defined by the short exact sequence
$$
\xymatrix{
0 \ar[r] &\mathcal O^{\mathrm{an}}(\mathcal F) \ar[r] & \mathcal R(\mathcal F) \ar[r] & \delta(\mathcal F) \ar[r] & 0.
}
$$
We will write $\mathcal O_{T}^{\mathrm{an}} := \Gamma(]T[, \mathcal O_{]T[})$ and call
$$
\mathcal R_{T} := (i_{T*}\mathcal O_{]T[})_{]\xi[} \quad \mathrm{and} \quad \delta_{T}:= \mathcal R_{T}/\mathcal O_{T}^{\mathrm{an}}
$$
the \emph{Robba ring} and the \emph{Dirac space}, respectively, of $T$.
\end{dfn}

\begin{prop}
If $x$ is a closed of point of $X_k$ with unramified lifting $a$, then $\mathcal R_{x}$ is (isomorphic to) the usual Robba ring $\mathcal R_{a}$ over $K(a)$.
In general, we have
$$
\mathcal R_{T}(\mathcal F) = \varinjlim \Gamma(]T[ \setminus ]T[_{\lambda}, \mathcal F)).
$$
\end{prop}

\textbf {Proof: }
As already mentioned, the subsets $X_{K}^{\mathrm{an}} \setminus ]Z[_{\lambda}$ form a cofinal family of affinoid neighborhoods of $]\xi[$ when $Z$ runs through the finite subsets of $X_k$ and $\lambda < 1$.
It follows that
$$
 \mathcal R_{T}(\mathcal F) = (i_{T*}\mathcal F_{|]T[})_{]\xi[} = \varinjlim \Gamma(X_{K}^{\mathrm{an}} \setminus ]Z[_{\lambda}, i_{T*}\mathcal F_{|]T[}) = \varinjlim \Gamma(]T[ \setminus ]T[_{\lambda}, \mathcal F)).
$$
In particular, we obtain
$$
\mathcal R_{x} = \varinjlim \Gamma(]x[ \setminus ]x[_{\lambda}, \mathcal O_{]x[})).
$$
Since $]x[ \simeq \mathbf D_{K(a)}(0, 1^-)$ and $]x[_{\lambda} \simeq \mathbf D_{K(a)}(0, \lambda^-)$, this is the (usual) Robba ring $\mathcal R_{a}$ over $K(a)$.
$\quad \Box$

If $T$ is a non empty finite closed subset of $X_{k}$, the adjunction map
$$
\mathcal O_{X_{K}^{\mathrm an}} \to i_{T*}i_{T}^{-1} \mathcal O_{X_{K}^{\mathrm an}}
$$
will induce on the stalks a canonical morphism $\mathcal O_{X_{K}^{\mathrm an},]\xi[} \to \mathcal R_{T}$.
In particular, if $U$ is a non empty affine open subset of $X_{k}$ with algebraic lifting $\mathrm{Spec}\; A$, there is a canonical morphism
$$
A^\dagger_{K} \to \mathcal O_{X_{K}^{\mathrm an},]\xi[} \to  \mathcal R_{T}.
$$

\begin{dfn} \label{robfib}
Let $T$ be a non empty finite closed subset of $X_{k}$.
If $E$ is a constructible module on $X_{K}^{\mathrm{an}}$, the \emph{Robba fiber} of $E$ at $T$ is
$$
R_{T}(E) := \mathcal R_{T} \otimes_{\mathcal O_{X_{K}^{\mathrm an},]\xi[}} E_{]\xi[}.
$$
Let $U$ be a non empty affine open subset of $X_{k}$ with algebraic lifting $ \mathrm{Spec}\; A$ and $M$ a finite $A^\dagger_{K}$-module.
Then, the \emph{Robba fiber} of $M$ at $T$ is
$$
R_{T}(M) := \mathcal R_{T} \otimes_{A^\dagger_{K}} M.
$$
\end{dfn}

Note that the Robba fiber is a free module (of finite rank) because $\mathcal O_{X_{K}^{\mathrm an},]\xi[}$ is a field.
Note also that if $\mathrm {Spec}\;B$ is some affine open subset of $\mathrm {Spec}\; A$, then
$$
\mathcal R_{T}(B^\dagger_{K} \otimes_{A^\dagger_{K}} M) = \mathcal R_{T}(M).
$$
Finally, if $i_{U}^{-1}E$ is coherent and $M=  \Gamma(]U[, i_{U}^{-1}E)$, then $\mathcal R_{T}(M) = \mathcal R_{T}(E)$.
Be careful however that $\mathcal R_{T}(E) \neq \mathcal R(E_{|]T[})$ in general (unless $E$ is coherent) but we have the following:

\begin{lem} \label{robbafib}
Let $U$ be a non empty affine open subset of $X_{k}$ with algebraic lifting $\mathrm{Spec}\; A$ and $T$ be a non empty finite closed subset of $U$.
Let $M$ be a finite $A^\dagger_{K}$-module and $M_{T}$ the restriction to $]T[$ of the corresponding $i_{U}^{-1} \mathcal O_{X_{K}^{\mathrm{an}}}$-module.
Then, the Robba fiber of $M$ at $T$ is identical to the Robba module of $M_{T}$:
$$
\mathcal R_{T}(M) = \mathcal R(M_{T}).
$$
\end{lem}

\textbf {Proof: }
Since we are working with right exact functors, we may assume that $M = A^\dagger_{K}$ in which case $M_{T} = \mathcal O_{]T[}$ and we fall back onto the definition.
$\quad \Box$

\begin{dfn}
Let $U$ be a non empty affine open subset of $X_{k}$ with closed complement $Z$ and $\mathrm{Spec}\; A$ an algebraic lifting of $U$.
Let $M$ be a finite $A^\dagger_{K}$-module.

An \emph{extension (at infinity)} of $M$ is a coherent module $\mathcal F$ on $]Z[$ together with a linear map $\mathcal R_{Z}(M) \to \mathcal R(\mathcal F)$.
It is said to be \emph{free} if $M$ is projective and $\mathcal F$ is locally free.

A {morphism of extensions} is pair made of a linear map $M' \to M$ and a linear map $\mathcal F' \to \mathcal F$ making commutative the obvious diagram.
\end{dfn}

In practice, we may just say that $\mathcal R_{Z}(M) \to \mathcal R(\mathcal F)$ `is' an extension, $M$ and $\mathcal F$ being understood as being part of the data.

\begin{dfn}
\begin{enumerate}
\item
Let $U$ be a non empty affine open subset of $X_{k}$ with closed complement $Z$ and $\mathrm{Spec}\; A$ an algebraic lifting of $U$.
Let $\mathcal R_{Z}(M) \to \mathcal R(\mathcal F)$ be an extension.
Let $T \subset U$ be any finite closed subset and $W = \mathrm {Spec}\; B$ some affine open subset of $V$ which is a lifting of $U \setminus T$.
Then, the extension of $N := B^\dagger_{K} \otimes_{A^\dagger_{K}} M$ defined by
$$
\mathcal R_{Z \cup T}(N) = \mathcal R_{Z \cup T}(M) =  \mathcal R_{Z}(M) \oplus \mathcal R_{T}(M) \to \mathcal R(\mathcal F) \oplus R(M_{T}).
$$
is called the \emph{shrinking} of $\mathcal R_{Z}(M) \to \mathcal R(\mathcal F)$ to $U \setminus T$.
\item Two extensions $\mathcal R_{Z}(M) \to \mathcal R(\mathcal F)$ on $U$ and $\mathcal R_{Z}(M') \to \mathcal R(\mathcal F')$ on $U'$ are said to be \emph{equivalent} if their shrinkings to $U \cap U'$ are isomorphic.
\end{enumerate}
\end{dfn}

\begin{prop} \label{smallthm}
The category of constructible (resp. constructible-free) modules on $X_{K}^{\mathrm{an}}$ is equivalent to the category of extensions (resp. free extensions) modulo equivalence.
\end{prop}

\textbf {Proof: }
If we are given a closed subset $Z$ of $X$ with closed complement $U$, we know that $]U[$ is a closed subset with open complement $]Z[$ and it is a general fact that the category of $\mathcal O_{X_{K}^{\mathrm{an}}}$modules $E$ on $X_{K}^{\mathrm{an}}$ is equivalent to the category of triples made of an $i_{U}^{-1} \mathcal O_{X_{K}^{\mathrm{an}}}$-module $E_{U}$, an $i_{Z}^{-1} \mathcal O_{X_{K}^{\mathrm{an}}}$-module $E_{Z}$ and an $i_{U}^{-1} \mathcal O_{X_{K}^{\mathrm{an}}}$-linear map $E_{U} \to i_{U}^{-1}i_{Z*}E_{Z}$.
This equivalence is described by the following morphism of short exact sequences:
$$
\xymatrix{
0 \ar[r] & i_{Z!}E_{Z} \ar[r] & i_{Z*}E_{Z} \ar[r] & i_{U*}i_{U}^{-1}i_{Z*}E_{Z} \ar[r] & 0 
\\
0 \ar[r]  & i_{Z!}E_{Z} \ar[r] \ar@{=}[u] & E \ar[r] \ar[u]  & i_{U*}E_{U} \ar[r] \ar[u] & 0 
.}
$$

We know from corollary \ref{locdec} that $E$ is constructible if and only if we can find such $U$ and $Z$ with both $E_{U}$ and $E_{Z}$ coherent.
Then, we simply set $\mathcal F = E_{Z}$ and $M = \Gamma(]U[, E_{U})$.
Now, if $i'_\xi : ]\xi[ \hookrightarrow ]U[$ denotes the inclusion map, we have
$$
i_U^{-1}i_{Z*}\mathcal F = i'_{\xi*}i'^{-1}_\xi i_U^{-1}i_{Z*}\mathcal F = i'_{\xi*}\mathcal R(\mathcal F)
$$
and it follows that
$$
\mathrm {Hom}_{i_{U}^{-1}\mathcal O_{X_{K}^{\mathrm{an}}}}(E_{U}, i_U^{-1}i_{Z*}E_{Z}) = \mathrm {Hom}_{i_{U}^{-1}\mathcal O_{X_{K}^{\mathrm{an}}}}(E_{U}, i'_{\xi*}\mathcal R(\mathcal F))
$$
$$
=  \mathrm {Hom}_{\mathcal O_{X_{K}^{\mathrm{an}},]\xi[}}(E_{U,]\xi[}, \mathcal R(\mathcal F)) = \mathrm {Hom}_{\mathcal R_{Z}}(\mathcal R_{Z}(M), \mathcal R(\mathcal F))
$$
because $E_{]U[,]\xi[} = \mathcal O_{X_{K}^{\mathrm{an}},]\xi[} \otimes_{A^\dagger_{K}} M$.
Thus, we see that morphisms $E_{U} \to i_{U}^{-1}i_{Z*}E_{Z}$ correspond bijectively to morphisms $\mathcal R_{Z}(M) \to \mathcal R(\mathcal F)$.

Of course, one easily checks that if we shrink $U$, the corresponding extension will be the shrinking of $\mathcal R_{Z}(M) \to \mathcal R(\mathcal F)$.
$\quad \Box$

In practice, a morphism $\mathcal R_{Z}(M) \to \mathcal R(\mathcal F)$ corresponds to a morphism $i_{U*}E_{U} \to i_{\xi*}\mathcal R(\mathcal F)$, and we can pull back the exact sequence
$$
0 \to i_{Z!}\mathcal F \to i_{Z*}\mathcal F \to i_{\xi*}R(\mathcal F) \to 0 
$$
in order to get
$$
0 \to i_{Z!}E_{Z} \to E \to i_{U*} E_{U} \to 0.
$$

\section{Specialization} \label{secsp}

The specialization map $\mathrm{sp} : X_{K}^{\mathrm{an}} \to X_{k}$ is not continuous for the usual topology of $X_{K}^{\mathrm{an}}$ (actually, it is anticontinuous) and it is therefore necessary to use the Tate topology in order to define direct image.

\begin{dfn}
If $E$ is a any module on $X_{K}^{\mathrm{an}}$, the \emph{specialization} of $E$ is
$$
\widetilde{ \mathrm{sp}}_{*} E := \mathrm{sp}_{*} E_{G}
$$
where $E_{G}$ denotes the module associated to $E$ for the Tate topology.
\end{dfn}

\begin{prop} \label{sperig}
The specialization functor $\widetilde{ \mathrm{sp}}_{*}$ is left exact.
Moreover, if $E$ is a any module on $X_{K}^{\mathrm{an}}$, then
$$
\mathrm R \widetilde{ \mathrm{sp}}_{*}E = \mathrm R \mathrm{sp}_{*} E_{G} = \mathrm R \mathrm{sp}_{*}E_{0}
$$
where $E_{G}$ (resp. $E_{0}$) denotes the module associated to $E$ for the Tate topology (on the corresponding rigid space).
\end{prop}

\textbf {Proof: }
The functor is left exact as composition of two left exact functors.
Moreover, since $E \mapsto E_{G}$ is even exact, we have $\mathrm R \widetilde{ \mathrm{sp}}_{*} E := \mathrm R \mathrm{sp}_{*} E_{G}$.
Finally, we have $\Gamma(]U[, E_{G}) = \Gamma(]U[^{\mathrm{rig}}, E_{0})$ if $U$ is any affine open subset of $X_{k}$.
$\quad \Box$

Recall now from \cite{Berthelot96}, 4.2.4, that if $T$ is a finite closed subset of $X_{k}$, then the sheaf $\mathcal O_{\widehat{X}}({}^\dagger T)$ may be defined as follows.
First of all, if $j : X_{k} \setminus Z \hookrightarrow X_{k}$ denotes the inclusion map, then $\mathcal O_{\widehat{X}}({}^\dagger T)$ is a subsheaf of $j_{*}j^{-1}\mathcal O_{\hat X}$.
Moreover, if $T$ is defined in $\mathcal U = \mathrm{Spf}\; \mathcal A$ by $f = 0 \mod \mathfrak m$, then
$$
\Gamma(U, \mathcal O_{\widehat{X}}({}^\dagger T)) := \varinjlim \mathcal A\{\pi/f^r\}
$$
so that
$$
\Gamma(U, \mathcal O_{\widehat{X}}({}^\dagger T)_{\mathbf Q}) := \mathcal A[1/f]^\dagger_{\mathbf Q}.
$$

The next result is a Berkovich analog of a theorem of Berthelot (proposition 4.3.2 of \cite{Berthelot96}, using proposition \ref{comprig}).

\begin{prop} \label{spope}
If $U$ is a non empty open subset of $X_k$ with closed complement $Z$, then
$$
\mathrm R\widetilde{\mathrm{sp}}_*i_{U*}i_U^{-1}\mathcal O_{X_K^{\mathrm{an}}} = \mathcal O_{\widehat{X}}({}^\dagger Z)_{\mathbf Q}.
$$
Moreover, the functor $\mathrm R\widetilde{\mathrm{sp}}_*i_{U*}$ induces an equivalence between coherent  $i_U^{-1}\mathcal O_{X_K^{\mathrm{an}}}$-modules and coherent  $\mathcal O_{\widehat{X}}({}^\dagger Z)_{\mathbf Q}$-modules.
\end{prop}

\textbf {Proof: }
If $E$ is a coherent $i_U^{-1}\mathcal O_{X_K^{\mathrm{an}}}$-module, we may always consider it as the restriction to $]U[$ of some coherent sheaf $\mathcal F$ defined on some neighborhood $V_{\lambda} = X_K^{\mathrm{an}} \setminus ]Z[_{\lambda}$ of $]U[$.
By definition, $\mathrm R^q\widetilde{\mathrm{sp}}_*i_{U*}E$ is the sheaf associated to
$$
U' \mapsto \mathrm H^q(]U'[, i_{U'}^{*}i_{U*}E)
$$
If we set $V'_{\lambda} = ]U'[ \cap V_{\lambda}$ and denote by $i_{\lambda} : V'_{\lambda} \hookrightarrow V_{\lambda}$ the inclusion map, it is a general topological result (proposition 2.5 of \cite{KashiwaraSchapira90} for example) that
$$
\mathrm H^q(]U'[, i_{U'}^{*}i_{U*}E) = \varinjlim \mathrm H^q(V'_{\lambda}, i_{\lambda}^*\mathcal F).
$$
Since $V'_{\lambda}$ is affine, we have $\mathrm H^q(]U'[, i_{U'}^{*}i_{U*}E)$ for $q > 0$ and therefore, $\mathrm R^q\widetilde{\mathrm{sp}}_*i_{U*}E_{U} = 0$ for $q > 0$.
Moreover, we have $V'_{\lambda} = ]U'[ \setminus ]Z'[_{\lambda}$ if $Z' = Z \cap U'$.
If we call $\mathcal U' = \mathrm{Spf}\mathcal A'$ the formal lifting of $U'$ and if $Z'$ is defined by $f = 0 \mod \mathfrak m$ in $\mathcal U'$, we have
$V'_{\lambda} = \mathrm{Spf}\mathcal A'_{K}\{\pi/f^r\}$ for $\lambda = |\pi|^{1/r}$ and therefore,
$$
\Gamma(]U'[, i_U^{-1}\mathcal O_{X_K^{\mathrm{an}}}) = \varinjlim \mathcal A'_{K}\{\pi/f^r\} = \Gamma(U,  \mathcal O_{\widehat{X}}({}^\dagger Z)_{\mathbf Q}).
$$
The last assertion then follows from theorem A and B for $\mathcal O_{\widehat{X}}({}^\dagger Z)_{\mathbf Q}$ (see proposition 4.3.2 of \cite{Berthelot96}).
$\quad \Box$

In practice, if  $\mathrm{Spec}\;A$ is an algebraic lifting of $U$, a coherent $i_U^{-1}\mathcal O_{X_K^{\mathrm{an}}}$-module $E$ corresponds to a finite $A^\dagger_{K}$-module M.
Now, let $U'$ be any affine open subset of $X$ with formal lifting $\mathcal U' = \mathrm{Spf}\; \mathcal A'$ and $f$ an equation for $Z$ in $\mathcal U'$.
Then, there is a canonical map $A^\dagger \to \mathcal A'[1/f]^\dagger$ 	and
$$
\Gamma(U', \widetilde{\mathrm{sp}}_*i_{U*} E) = \mathcal A'[1/f]^\dagger_{\mathbf Q} \otimes_{A^\dagger_{K}} M.
$$

We now study the case of a finite closed subset $Z \subset X_{k}$.
We will write
$$
\mathcal H^\dagger_Z := \mathcal O_{\widehat{X}}({}^\dagger Z)/\mathcal O_{\widehat{X}}.
$$
For those who now the theory of arithmetic $\mathcal D$-modules, we recall that
$$
\mathrm R\underline \Gamma^\dagger_{Z} \mathcal O_{\hat X\mathbf Q} \simeq \mathcal H^\dagger_{Z\mathbf Q} [-1].
$$
This cohomology with support is closely related to the notion of Dirac space introduced in definition \ref{robdir} as we shall see right now.

\begin{lem}\label{spclo}
If $Z$ is a finite closed subset of $X_{k}$, we have
$$
\mathcal H^\dagger_{Z\mathbf Q} = \bigoplus_{x \in Z}i_{x*} \delta_{x}
$$
where $i_{x} : \{x\} \hookrightarrow \hat X$ denotes the inclusion map.
\end{lem}

\textbf {Proof: }
Since $\mathcal O_{\widehat{X}}({}^\dagger Z)$ and $\mathcal O_{\widehat{X}}$ coincide outside $Z$, it is sufficient to show that there is an isomorphism on the stalks $\mathcal H^\dagger_{Z\mathbf Q,x} \simeq \delta_{x}$ when $x \in Z$.
We will now use the fact that if $U$ is an affine neighborhood of $x$, then $]U[ \setminus ]x[_{\lambda}$ and $]x[$ form an open covering of $]U[$ with intersection $]x[ \setminus ]x[_{\lambda}$ and that this covering is acyclic for $\mathcal O_{]U[}$.
We call $\mathcal U = \mathrm{Spf}\; \mathcal A$ the formal lifting of $U$ and assume that $x$ is defined in $\mathcal U$ by an equation $f = 0 \mod \mathfrak m$.
We will also denote by $\mathcal R_{x,\lambda} := \Gamma(]x[ \setminus ]x[_{\lambda}, \mathcal O_{]x[})$ and $\delta_{x,\lambda} = \mathcal R_{x,\lambda}/\mathcal O_{x}^{\mathrm{an}}$ with $\mathcal O_{x}^{\mathrm{an}} = \Gamma(]x[, \mathcal O_{]x[})$ as before.
Then, the above remark implies that, in the following morphism of exact sequences, the last map is an isomorphism
$$
\xymatrix{
0 \ar[r] & \mathcal A_{\mathbf Q} \ar[r] \ar[d]& \mathcal A\{\pi/f^r\}_{\mathbf Q} \ar[r] \ar[d] & (\mathcal A\{\pi/f^r\}/\mathcal A)_{\mathbf Q} \ar[r] \ar[d]^\simeq & 0 
\\
0 \ar[r] & \mathcal O^{\mathrm{an}}_{x} \ar[r] & \mathcal R_{x,\lambda} \ar[r] & \delta_{x,\lambda} \ar[r] & 0
}
$$
with $\lambda = |\pi|^{1/r}$.
Taking limit on $\lambda < 1$ (and all $U \ni x$), we obtain the expected isomorphism $\mathcal H^\dagger_{Z\mathbf Q,x} = \mathcal \delta_{x}$.
$\quad \Box$

\begin{lem}\label{spclo2}
If $Z$ is a non empty finite closed subset of $X_{k}$, then
$$
\mathrm R\widetilde{\mathrm{sp}}_*i_{Z!}\mathcal O_{]Z[} \simeq \mathcal H^\dagger_{Z\mathbf Q} [-1] \quad ( = \mathrm R\underline \Gamma^\dagger_{Z} \mathcal O_{\hat X\mathbf Q}).
$$
\end{lem}

\textbf {Proof: }
If $U$ denotes the complement of $Z$ in $X_{k}$, there is an exact sequence
$$
0 \to i_{Z!}i_{Z}^{-1}\mathcal O_{X_K^{\mathrm{an}}} \to \mathcal O_{X_K^{\mathrm{an}}} \to i_{U*} i_{U}^{-1}\mathcal O_{X_K^{\mathrm{an}}} \to 0 
$$
from which we derive a triangle
$$
\mathrm R\widetilde{\mathrm{sp}}_* i_{Z!}\mathcal O_{]Z[} \to \mathrm R\widetilde{\mathrm{sp}}_* \mathcal O_{X_K^{\mathrm{an}}} \to \mathrm R\widetilde{\mathrm{sp}}_* i_{U*} i_{U}^{-1}\mathcal O_{X_K^{\mathrm{an}}} \to
$$
and then, using proposition 3.2, we obtain
$$
\mathrm R\widetilde{\mathrm{sp}}_* i_{Z!}\mathcal O_{]Z[} \simeq [\mathcal O_{\widehat{X}\mathbf Q} \to \mathcal O_{\widehat{X}}({}^\dagger Z)_{\mathbf Q}] \simeq \mathcal H^\dagger_{Z\mathbf Q} [-1].
\quad \Box$$

Recall that we can always lift a finite closed subset $Z$ to a smooth relative divisor $D \subset X$: if $Z = \{x_{1}, \ldots, x_{n}\}$, we choose an unramified lifting $a_{i} \in X_{K}$ for each $i$ and let $D =\mathrm{Spec}\; \prod \mathcal V(a_{i})$.
Note then that a (coherent) $\mathcal O_{D_{K}}$-module $H$ is simply a finite collection of (finite dimensional) $K(a)$-vector spaces $H_{a}$ for $a \in D_{K}$.
Note also that the canonical inclusion $D_{K} = D_{K}^{\mathrm{an}}= \hat D_{K} \hookrightarrow ]Z[$ has a unique retraction contracting each disc onto its ``center''.
It also follows from lemma \ref{spclo} that $\mathcal H^\dagger_{Z\mathbf Q}$ is an $\mathcal O_{D_{K}}$-module.

\begin{prop} \label{pointwis}
Let $D \subset X$ be a smooth divisor with reduction $Z$.
If $E$ is a locally free module on $]Z[$, there exists a coherent $\mathcal O_{D_{K}}$-module $H$ such that $E \simeq \mathcal O_{]Z[} \otimes_{\mathcal O_{D_{K}}} H$ and we have
$$
\mathrm R\widetilde{\mathrm{sp}}_*i_{Z!}E \simeq (\mathcal H^\dagger_{Z\mathbf Q} \otimes_{\mathcal O_{D_{K}}} H)[-1]
$$
\end{prop}

\textbf {Proof: }
The first assertion follows from the remark before and the second one then follows from lemma \ref{spclo2} since $\mathrm R\widetilde{\mathrm{sp}}_*$ is additive.
$\quad \Box$

\begin{dfn} \label{perv} A \emph{perverse sheaf} on $\hat X$ is a complex of $\mathcal O_{\hat X\mathbf Q}$-modules $\mathcal E$ such that $\mathcal H^{0}(\mathcal E)$ is $\mathcal O_{\hat X\mathbf Q}$-flat, $\mathcal H^{1}(\mathcal E)$ has finite support and $\mathcal H^{i}(\mathcal E) = 0$ otherwise.
\end{dfn}

It is sometimes convenient to split this definition in two (in order to see perverse sheaves as the heart of a $t$-structure).
We can denote by $D_{\geq 0}(\hat X)_{\mathbf Q}$ the category of bounded complexes of $\mathcal O_{\hat X\mathbf Q}$-modules $\mathcal E$ where $\mathcal H^{0}(\mathcal E)$ is $\mathcal O_{\hat X\mathbf Q}$-flat and $\mathcal H^{i}(\mathcal E) = 0$ for $i < 0$.
Next, we denote by $D_{\leq 0}(\hat X) _{\mathbf Q} $ the category of bounded complexes of $\mathcal O_{\hat X\mathbf Q}$-modules $\mathcal E$ where $\mathcal H^{1}(\mathcal E)$ has finite support and $\mathcal H^{i}(\mathcal E) = 0$ for $i > 1$.
Then the category of perverse sheaves is $D_{\geq 0}(\hat X)_{\mathbf Q} \cap D_{\leq 0}(\hat X) _{\mathbf Q} $.

\begin{prop} \label{descdir}
If $E$ is a constructible-free module on $X_{K}^{\mathrm{an}}$, then $\mathrm R\widetilde{\mathrm{sp}}_* E$ is a perverse sheaf on $\hat X$.
More precisely, we have
$$
\mathrm R\widetilde{\mathrm{sp}}_* E \simeq [\mathcal E \to \mathcal H^\dagger_{Z} \otimes_{\mathcal O_{D_{K}}} H]
$$
where $Z \subset X_{k}$ is a finite closed subset, $\mathcal E$ is a coherent locally free $\mathcal O_{\widehat{X}}({}^\dagger Z)_{\mathbf Q}$-module, $D \subset X$ is a smooth divisor with reduction $Z$,  and $H$ is a coherent $\mathcal O_{D_{K}}$-module.
\end{prop}

\textbf {Proof: }
We showed in corollary \ref{excons} that there exists an exact sequence
$$
0 \to i_{Z!}E_{Z} \to E \to i_{U*} E_{U} \to 0 
$$
where $Z \subset X_{k}$ is a finite closed subset with affine open complement $U$, $E_{U}$ is a coherent locally free $i_{U}^{-1} \mathcal O_{X_{K}^{\mathrm{an}}}$-module and $E_{Z}$ is a locally free $\mathcal O_{]Z[}$-module.

There exists a smooth divisor $D \subset X$ with reduction $Z$ and a coherent $\mathcal O_{D_{K}}$-module $H$ such that $E \simeq \mathcal O_{]Z[} \otimes_{\mathcal O_{D_{K}}} H$.
We saw in proposition \ref{pointwis} that
$$
\mathrm R\widetilde{\mathrm{sp}}_*i_{Z!}E \simeq (\mathcal H^\dagger_{Z} \otimes_{\mathcal O_{D_{K}}} H)[-1].
$$
On the other hand, we know from proposition \ref{spope} that
$$
\mathrm R\widetilde{\mathrm{sp}}_*i_{U*} E_{U} = \mathcal E
$$
where $\mathcal E$ is a coherent locally free $\mathcal O_{\widehat{X}}({}^\dagger Z)_{\mathbf Q}$-module.
From the above exact sequence, we obtain the following exact triangle
$$
( \mathcal H^\dagger_{Z} \otimes_{K_{Z}} H) [-1] \to \mathrm R\widetilde{\mathrm{sp}}_* E \to \mathcal E \to
$$
and we are done.
$\quad \Box$

In practice, a constructible module $E$ is given by some extension $\mathcal R(M) \to R(\mathcal F)$ where $M$ is a finite $A^\dagger_{K}$ module ($\mathrm{Spec}\;A$ being an algebraic lifting of the open complement $U$ of $Z$ in $X_{k}$) and $\mathcal F$ is a coherent module on $]Z[$.
Now, let $U'$ be any affine open subset of $X$ with formal lifting $\mathcal U' = \mathrm{Spf}\; \mathcal A'$ and $f$ an equation for $Z$ in $\mathcal U'$.
Let $Z' := U' \cap Z$ and $\mathcal F' := i_{Z'}^{-1}E$.
There are canonical maps $A^\dagger \to \mathcal A'[1/f]^\dagger$ as well as $\mathcal A'[1/f]^\dagger_{\mathbf Q} \to R_{x}$ when $x \in U'$.
We can compose the map
$$
\mathcal A'[1/f]^\dagger_{\mathbf Q} \otimes_{A^\dagger_{K}} M \to \mathcal R(\mathcal F')
$$
obtained by scalars extensions from $M \to \mathcal R(M) \to \mathcal R(\mathcal F) \to  \mathcal R(\mathcal F')$ with the projection $\mathcal R(\mathcal F') \to \delta(\mathcal F')$.
When $E$ is constructible-free (the other cases are not interesting for us), we get
$$
\Gamma(U', \mathrm R\widetilde{\mathrm{sp}}_* E) = \left[\mathcal A'[1/f]^\dagger_{\mathbf Q} \otimes_{A^\dagger_{K}} M \to \delta(\mathcal F')\right].
$$

Again, for latter use, we show that specialization commutes with Frobenius.

\begin{prop}
Assume that $F: X \to X$ is a lifting of the absolute Frobenius of $X_{k}$.
If $E$ is an $\mathcal O_{X_{K}^{\mathrm{an}}}$-module, there is a natural isomorphism $\hat F^*\mathrm R\widetilde{\mathrm{sp}}_* E \simeq \mathrm R\widetilde{\mathrm{sp}}_* F_{K}^{\mathrm{an}*}E$.
\end{prop}

\textbf {Proof: }
We have to be a little careful because $\widetilde{\mathrm{sp}}_*$ is not a direct image.
The method is to start directly with an injective resolution $I$ of $E_{G}$.
It is then sufficient to check that $\hat F^*\widetilde{\mathrm{sp}}_* E \simeq \widetilde{\mathrm{sp}}_* F_{K}^{\mathrm{an}*}E$ and this easily follows from the fact that $F$ is finite and flat.
$\quad \Box$

\section{Connections on constructible modules}

We start with the very general following definition.

\begin{dfn}
If $Y \subset X_{k}$ is any locally closed subset and $i_{Y} : ]Y[ \hookrightarrow X_{K}^{\mathrm{an}}$ denotes the inclusion map, a \emph{connection} on an $i_{Y}^{-1}\mathcal O_{X_{K}^{\mathrm{an}}}$-module $E$ is a $K$-linear map
$$
\nabla : E \to E \otimes_{i_{Y^{-1}}\mathcal O_{X_{K}^{\mathrm{an}}}} i_{Y}^{-1}\Omega^1_{X_{K}^{\mathrm{an}}}
$$
satisfying the Leibnitz rule.
A \emph{$\nabla$-module} on $]Y[$ is an $i_{Y}^{-1}\mathcal O_{X_{K}^{\mathrm{an}}}$-module $E$ endowed with a connection.
A \emph{horizontal map} between $\nabla$-modules on $Y$ is a $i_{Y}^{-1}\mathcal O_{X_{K}^{\mathrm{an}}}$-linear map that commutes with the connection.
\end{dfn}

Of course, if $Y' \subset Y$, any $\nabla$-module $E$ on $]Y[$ restricts to a $\nabla$-module on $]Y'[$.
Note also that if $Z \subset X_{k}$ is a closed subset, meaning either a finite subset or $X_{k}$ itself, then $]Z[$ is open in $X_{K}^{\mathrm{an}}$ and we fall back onto the usual notion of connection in analytic geometry.
Finally, note - and this is important - that there is \emph{no} finiteness condition at this point.

\begin{prop}
A constructible $\nabla$-module $E$ on $X_{K}^{\mathrm{an}}$ is constructible-free.
\end{prop}

\textbf {Proof: }
Since a coherent module with a connection on a disc is necessarily free when the valuation is discrete, we can apply corollary \ref{freeness}.
$\quad \Box$

It is very convenient to be able to use the description of connections in terms of stratifications.
We recall the definition of the first infinitesimal neighborhood $P$ of the diagonal.
If $X$ is defined by some ideal $\mathcal I$ into $X \times X$, then $P$ is the closed subscheme of $X$ defined by $\mathcal I^2$.
Recall that $X$ and $P$ have same underlying space and that there is a short exact sequence
$$
0 \to  \Omega^1_{X} \to \mathcal O_{P} \to  \mathcal O_{X} \to 0.
$$
Actually, we will only need the analogous construction on $X_{K}^{\mathrm{an}}$ which can also be deduced from this one by functoriality.
We will consider the maps $q_{1}, q_{2} = P_{K}^{\mathrm{an}} \to X_{K}^{\mathrm{an}}$ and $\Delta : X_{K}^{\mathrm{an}} \hookrightarrow P_{K}^{\mathrm{an}}$ induced by the projections and the diagonal embedding.

\begin{dfn}
If $Y \subset X_{k}$ is a locally closed subset, a \emph{1-stratification} on an $i_{Y}^{-1}\mathcal O_{X_{K}^{\mathrm{an}}}$-module $E$ is an isomorphism
$$
\epsilon : i_{Y}^{-1}q_{2}^*i_{Y*}E \simeq i_{Y}^{-1}q_{1}^*i_{Y*} E
$$
such that $(i_{Y}^{-1}\Delta^*i_{Y*})(\epsilon)$ is the identity of $E$.
A morphism of such is an $i_{Y}^{-1}\mathcal O_{X_{K}^{\mathrm{an}}}$-linear map that commutes with the 1-stratifications.
\end{dfn}

Alternatively, we can write this isomorphism as
$$
i_{Y}^{-1}\mathcal O_{P_{K}^{\mathrm{an}}} \otimes_{i_{Y}^{-1}\mathcal O_{X_{K}^{\mathrm{an}}}} E \simeq E \otimes_{i_{Y}^{-1}\mathcal O_{X_{K}^{\mathrm{an}}}} i_{Y}^{-1}\mathcal O_{P_{K}^{\mathrm{an}}}.
$$
The induced map $E \to E \otimes_{i_{Y}^{-1}\mathcal O_{X_{K}^{\mathrm{an}}}} i_{Y}^{-1}\mathcal O_{P_{K}^{\mathrm{an}}}$ coincides with the inclusion mod $i_{Y}^{-1} \mathcal I_{K}^{\mathrm{an}}$ and the difference is a map
$$
\nabla : E \to E \otimes_{i_{Y^{-1}}\mathcal O_{X_{K}^{\mathrm{an}}}} i_{Y}^{-1}\Omega^1_{X_{K}^{\mathrm{an}}}.
$$

\begin{prop}
This construction establishes an equivalence between the category of  $i_{Y}^{-1}\mathcal O_{X_{K}^{\mathrm{an}}}$-modules endowed with a 1-stratification and the category of $\nabla$-modules on $Y$.
\end{prop}

\textbf {Proof: }
Standard.
$\quad \Box$

Note that the category of $\nabla$-modules on $Y$ is also equivalent to the category of left $i_{Y}^{-1}\mathcal D_{X_{K}^{\mathrm{an}}}$-modules if $\mathcal D_{X_{K}^{\mathrm{an}}}$ denotes the sheaf of (algebraic) differential operators on $X_{K}^{\mathrm{an}}$.

Let $U \subset X_{k}$ be an affine open subset and $\mathrm{Spec}\; A \subset X$ an algebraic lifting.
We will now recall the description of a connection on an $A^\dagger_{K}$-module in term of stratification as we just did above for $i_{Y}^{-1}\mathcal O_{X_{K}^{\mathrm{an}}}$-modules.
If we set $P_{A} := (A \otimes A)/I^2$, where $I$ is the ideal of multiplication $A \otimes A \to A$, then a {1-stratification} on an $A^\dagger_{K}$-module $M$ a $P_{AK}^{\dagger}$-linear isomorphism
$$
P_{AK}^{\dagger} \otimes_{A^\dagger_{K}} M \simeq M \otimes_{A^\dagger_{K}} P_{AK}^{\dagger} 
$$
Again, this is equivalent to a connection on $M$.

\begin{prop}
Let $U \subset X_{k}$ be an affine open subset of $X_{k}$ with algebraic lifting $\mathrm{Spec}\; A$.
Then, the functor $\Gamma(]U[, -)$ induces an equivalence between coherent (necessarily locally free) $\nabla$-modules on $U$ and finite (necessarily projective) $\nabla$-$A^\dagger_{K}$-modules.
\end{prop}

\textbf {Proof: }
We showed in proposition \ref{dageq} that  $\Gamma(]U[, -)$ induces an equivalence between coherent  modules $E$ on $U$ and finite $A^\dagger_{K}$-modules $M$.
It follows that there is an equivalence between 1-stratifications on $E$ and 1-stratifications on $M$ because, obviously, with the above notations, we have
$$
\Gamma(]U[, i_{U}^{-1}\mathcal O_{P_{K}^{\mathrm{an}}}) = P_{AK}^{\dagger}.
\quad \Box$$

\begin{prop}
If $Y$ is a locally closed subset of $X_{k}$ and $E$ is a $\nabla$-module on $]Y[$, then $i_{Y*}E$ has a natural connection.
If $Z$ is a closed subset of $X_{k}$ and $E$ is a $\nabla$-module on $]Z[$, then $i_{Z!}E$ has a natural connection.
These functors are fully faithful.
\end{prop}

\textbf {Proof: }
We have for $j = 1, 2$,
$$
q_{j}^*i_{Y*} E \simeq i_{Y*}i_{Y}^{-1}q_{j}^*i_{Y*}E.
$$
Therefore, the 1-stratification extends canonically.
The partial inverse is as usual induced by $i_{Y}^{-1}$.
The proof follows the same lines for $i_{Z!}$ since we will have
$$
q_{j}^*i_{Z!} E \simeq i_{Z!}i_{Z}^{-1}q_{j}^*i_{Z!}E.
\quad \Box$$

\begin{prop} \label{exactseq}
A $\nabla$-module $E$ on $X_{K}^{\mathrm{an}}$ is constructible if and only if there exists an exact sequence
$$
0 \to i_{Z!}E_{Z} \to E \to i_{U*} E_{U} \to 0 
$$
where $U$ is a non-empty affine open subset of $X_{k}$ with closed complement $Z$, $E_{U}$ is a coherent (necessarily locally free) $\nabla$-module on $]U[$ and $E_{Z}$ is a coherent (necessarily locally free) $\nabla$-module on $]Z[$.
We can even assume that $E_{U}$ is free.
\end{prop}

\textbf {Proof: }
Using 1-stratifications, it follows from proposition \ref{excons}. 
$\quad \Box$

Recall that we defined above, when $T \subset X_{k}$ is a non empty finite closed subset and $\mathcal F$ is a coherent $\mathcal O_{]T[}$-module, the Robba module of $\mathcal F$ as
$$
\mathcal R(\mathcal F) := (i_{T*}\mathcal F)_{]\xi[}
$$
where $i_{T} : ]T[ \hookrightarrow X_{K}^{\mathrm{an}}$ is the inclusion map and $]\xi[$ is the ``generic point'' of $X_{K}^{\mathrm{an}}$.
In particular, we will consider the \emph{Robba ring} $\mathcal R_{T} := \mathcal R(\mathcal O_{T})$.
If for each $x \in T$, we choose an unramified lifting $a$ of $x $, then $\mathcal R_{T}$ is the direct sum of the usual Robba rings over $K(a)$.
We will denote by $\Omega^1_{R_{T}}$ the module of finite differentials over $\mathcal R_{T}$.

\begin{lem}
If $T$ is a non empty finite closed subset of $X_{k}$, then $\Omega^1_{\mathcal R_{T}}$ is canonically isomorphic to $\mathcal R(\Omega^1_{]T[})$.
Moreover, any connection on a coherent $\mathcal O_{]T[}$-module $\mathcal F$ induces an connection $\mathcal R(\mathcal F)$.
\end{lem}

\textbf {Proof: }
We may assume that $T$ is reduced to one rational point $x$.
Then, we know that $]x[$ is a disc with some parameter $t$, that $\mathcal R_{x}$ is the usual Robba ring and we have
$$
\mathcal R(\Omega^1_{]x[}) = \mathcal R(\mathcal O_{]x[} \mathrm dt) = \mathcal R_{x} \mathrm dt = \Omega^1_{\mathcal R_{x}}
$$
Also, we have
$$
\mathcal R(i_{x}^{-1}\mathcal O_{P_{K}^{\mathrm{an}}}) \simeq \mathcal R_{x}[\tau]/\tau^2
$$
where $\tau = p_{2}^{-1}(t) - p_{1}^{-1}(t)$.
Therefore a 1-stratification on an $\mathcal O_{]x[}$-module $\mathcal F$ will induce an isomorphism
$$
\mathcal R_{x}[\tau]/\tau^2 \otimes_{\mathcal R_{x}} \mathcal R(\mathcal F) \simeq \mathcal R(\mathcal F) \otimes_{\mathcal R_{x}} \mathcal R_{x}[\tau]/\tau^2.
$$
This is a 1-stratification on $\mathcal R(\mathcal F)$ that corresponds to a connection.
$\quad \Box$

Recall that if $T \subset X_{k}$ is a non empty finite closed subset and $U \subset X_{k}$ is a non empty affine open subset with algebraic lifting $\mathrm{Spec}\; A$, there is a canonical morphism $A^\dagger_{K} \to \mathcal R_{T}$ and that the \emph{Robba fiber} of an $A^\dagger_{K}$-module $M$ at $T$ is
$$
R_{T}(M) := \mathcal R_{T} \otimes_{A^\dagger_{K}} M.
$$
Note that if $M$ is endowed with a connection, then $R_{T}(M)$ inherits automatically a connection.

\begin{lem}
Let $U$ be a non empty affine open subset of $X_{k}$ with algebraic lifting $\mathrm{Spec}\; A$ and $T$ be a finite closed subset of $U$.
Let $M$ be a finite $\nabla$-$A^\dagger_{K}$-module and $M_{T}$ the restriction to $]T[$ of the corresponding $i_{U}^{-1} \mathcal O_{X_{K}^{\mathrm{an}}}$-module.
Then, there is a horizontal isomorphism
$$
\mathcal R_{T}(M) = \mathcal R(M_{T}).
$$
\end{lem}

\textbf {Proof: }
Deduced from lemma \ref{robbafib} using 1-stratifications.
$\quad \Box$

\begin{dfn} \label{nablext}
Let $U$ be a non empty affine open subset of $X_{k}$ with closed complement $Z$ and $\mathrm{Spec}\; A$ an algebraic lifting of $U$.
Let $M$ be a finite $\nabla$-$A^\dagger_{K}$-module.

A \emph{$\nabla$-extension (at infinity)} of $M$ is a coherent $\nabla$-module $\mathcal F$ on $]Z[$ together with a horizontal map $\mathcal R_{Z}(M) \to \mathcal R(\mathcal F)$.

A {morphism of $\nabla$-extensions} is a morphism of extensions given by horizontal maps.

Two $\nabla$-extensions $\mathcal R_{Z}(M) \to \mathcal R(\mathcal F)$ on $U$ and $\mathcal R_{Z}(M') \to \mathcal R(\mathcal F)$ on $U'$ are \emph{equivalent} if their shrinkings to $U \cap U'$ are isomorphic as $\nabla$-extensions.
\end{dfn}

\begin{prop}\label{equinab}
The category of constructible $\nabla$-modules on $X_{K}^{\mathrm{an}}$ is equivalent to the category of $\nabla$-extensions modulo equivalence.
\end{prop}

\textbf {Proof: }
Using proposition \ref{smallthm}, it is essentially sufficient to notice that if $E$ is a constructible module, then $q_{2}^{*}E$ too is constructible.
Actually, if $E$ is given by some extension $\mathcal R_{Z}(M) \to \mathcal R(\mathcal F)$, then $q_{2}^{*}E$ will be given by
$$
\mathcal R_{Z}(P_{AK}^\dagger \otimes_{A_{K}^\dagger }M) \to \mathcal R(i_{Z}^{-1}q_{2}^*\mathcal F)
$$
And the analogous results holds for $q_{1}^{*}E$.
Therefore, a 1-stratification on $E$ is equivalent to a 1-stratification on $M$ and a compatible 1-stratification on $\mathcal F$. 
$\quad \Box$

We denote by $\mathcal D_{\hat X}$ the sheaf of (algebraic) differential operators on $\hat X$.
We call a $\mathcal D_{\hat X\mathbf Q}$-module \emph{perverse} if the underlying $\mathcal O_{\hat X\mathbf Q}$ is perverse (see definition \ref{perv}).

\begin{prop}\label{nablam}
$\mathrm R\widetilde{\mathrm{sp}}_*$ induces a functor from the category of  constructible $\nabla$-modules on $X_{K}^{\mathrm{an}}$ to the category of perverse $\mathcal D_{\hat X\mathbf Q}$-modules.
\end{prop}

\textbf {Proof: }
Is is essentially sufficient to recall the construction of $\mathrm R\widetilde{\mathrm{sp}}_* E$ from the corresponding extension $\mathcal R(M) \to R(\mathcal F)$.
More precisely, let $U'$ be an affine open subset of $X$ with formal lifting $\mathcal U' = \mathrm{Spf}\; \mathcal A'$ and $f$ an equation for $Z$ in $\mathcal U'$.
Let $Z' := U' \cap Z$ and $\mathcal F' := i_{Z'}^{-1}\mathcal F$.
Then, we have
$$
\Gamma(U', \mathrm R\widetilde{\mathrm{sp}}_* E) = \left[\mathcal A'[1/f]^\dagger_{\mathbf Q} \otimes_{A^\dagger_{K}} M  \to \delta(\mathcal F')\right].
$$
where the non-trivial map in this complex is obtained by scalar extension from the composite $M \to \mathcal R(M) \to \mathcal R(\mathcal F)$, restriction of the image to the direct factor $\mathcal R(\mathcal F')$, and then, projection onto $\delta(\mathcal F')$.
All those maps are horizontal.
$\quad \Box$
 
We also want to state another result that we will need later on.
If $Z$ is a finite closed subset of $\hat X$, we set
$$
\mathcal D_{\hat X}({}^\dagger Z)_{\mathbf Q} = \mathcal O_{\hat X}({}^\dagger Z)_{\mathbf Q} \otimes_{\mathcal O_{\hat X}} \mathcal D_{\hat X}.
$$

\begin{prop} \label{weakequi}
If $U$ is a non empty affine open subset of $X_{k}$ with closed complement $Z$, then $\widetilde{\mathrm{sp}}_*i_{U*}$ establishes an equivalence between the category of coherent $\nabla$-modules on $]U[$ and $\mathcal D_{\hat X}({}^\dagger Z)_{\mathbf Q}$-modules that are coherent on $\mathcal O_{\hat X}({}^\dagger Z)_{\mathbf Q}$.
\end{prop}

\textbf {Proof: }
Use 1-stratifications again and proposition \ref{spope}.
$\quad \Box$

\begin{prop} 
Assume that $F: X \to X$ is a lifting of the absolute Frobenius of $X_{k}$.
\begin{enumerate}
\item If $Y$ is a locally closed subset of $X_{k}$ and $E$ is a $\nabla$-module on $]Y[$, then $i_{Y}^{-1}F_{K}^{\mathrm{an}*}i_{Y*}E$ is also a $\nabla$-module on $]Y[$.
\item If $E$ is a $\nabla$-module on $X_{K}^{\mathrm{an}}$, then the isomorphism
$$
\hat F^*\mathrm R\widetilde{\mathrm{sp}}_* E \simeq \mathrm R\widetilde{\mathrm{sp}}_* F_{K}^{\mathrm{an}*}E
$$
is horizontal.
\end{enumerate}
\end{prop}

\textbf {Proof: }
If we denote by $F_{P}$ the map induced on $P$ by $F \times F$, then for $j = 1, 2$, we have
$$
i_{Y}^{-1}q_{j}^*i_{Y*}i_{Y}^{-1} F_{K}^{\mathrm{an}*} i_{Y*} = i_{Y}^{-1}F_{PK}^{\mathrm{an}*}i_{Y*}i_{Y}^{-1}q_{j}^*i_{Y*}
$$
and the first assertion follows.
The second assertion is obtained by linearity by using 1-stratifications as usual.
$\quad \Box$

\begin{dfn}
Let $F: X \to X$ be an algebraic lifting of the absolute Frobenius of $X_{k}$ and $Y$ a locally closed subset of $X_{k}$.
An \emph{$F$-$\nabla$-module} on $]Y[$ is a $\nabla$-module $E$, endowed with an isomorphism $\phi : i_{Y}^{-1}F_{K}^{\mathrm{an}*}i_{Y*}E \simeq E$.
\end{dfn}

We let the reader extend all the results of this section to $F$-$\nabla$-modules.
This is straightforwards since $F$ is finite flat.

\section{Overconvergent connections}
 
Now, we embed $X$ diagonally into $X \times X$.
If $Y \subset X_{k}$ is any subset, we will denote by $]]Y[[$ the tube of  $Y$ into $X_{K}^{\mathrm{an}} \times X_{K}^{\mathrm{an}}$.
We call
$$
p_{1}, p_{2} : ]]X_{k}[[ \to X_{K}^{\mathrm{an}}
$$
the maps induced by the projections.
We will still denote by $i_{Y} : ]]Y[[ \hookrightarrow ]]X_{k}[[$ the embedding, hoping that this will not create any confusion.
Note that
$$
]]Y[[ = p_{j}^{-1}(]Y[) \subset ]]X_{k}[[ \quad \mathrm{for} \quad j = 1, 2.
$$
Sometimes, even if we will try to avoid it, we may also denote by $p_{1}, p_{2} : ]]Y[[ \to ]Y[$ the maps induced by the projections (at least when $Y$ is a closed subset).

Locally, the geometry of $]]X_{k}[[$ is not too bad (this is the strong fibration theorem of Berthelot) as we can see right now.
Since $X$ is smooth of relative dimension 1, there exists locally an \'etale map $t$ to the affine line $\mathbf A^1_{\mathcal V}$.
We will then say that $t$ is a \emph{local parameter} on $X$.
Assume that $t$ is defined on some open subset $V \subset X$ with reduction $U \subset X_{k}$.
Then $\tau := p_{2}^*(t) - p_{1}^*(t)$ defines an \'etale map $V \times V \to V \times \mathbf A^1_{\mathcal V}$.
More precisely, there is a commutative diagram
$$
\xymatrix{V \times V \ar[rr] \ar[rd]^{p_{1}} && V \times \mathbf A^1_{\mathcal V} \ar[ld]_{p}
\\
& V \ar@/^1pc/[lu]^{\Delta} \ar@/_1pc/[ru]_{0}&
}
$$
As in proposition \ref{pointdisc}, it follows from lemma 4.4 of \cite{Berkovich99} that this map induces an isomorphism $]]U[[ \simeq ]U[ \times \mathbf D(0, 1^-)$.
Moreover, the morphism
$$
V_{K}^{\mathrm{an}} \times V_{K}^{\mathrm{an}} \to V_{K}^{\mathrm{an}} \times \mathbf A_{K}^{1,\mathrm{an}}
$$
is \'etale (formally \'etale and boundaryless).
Thus, Proposition 4.3.4 of \cite{Berkovich93} implies that it induces an isomorphism between a neighborhood $V'$ of $]]U[[$ in $]]X_{k}[[$ and a neighborhood $V''$ of $]U[ \times \mathbf D(0, 1^-)$ into $X_{K}^{\mathrm{an}} \times \mathbf D(0, 1^-)$.
If we denote by $Z$ the closed complement of $U$ and let $V_{\lambda} = X_{K}^{\mathrm{an}} \setminus ]Z[_{\lambda}$ as usual, we may assume that $V'' = \cup_{\lambda, \eta} (V_{\lambda} \times \mathbf D(0, \eta^+))$ with $\lambda, \eta \to 1$.
We can summarize the situation in the following diagram
$$
\xymatrix{
V_{K}^{\mathrm{an}} \times V_{K}^{\mathrm{an}} \ar[r] & V_{K}^{\mathrm{an}} \times \mathbf A_{K}^{1,\mathrm{an}}
\\
V' \ar[r]^-\simeq \ar@{^{(}->}[u] & \cup_{\lambda, \eta} (V_{\lambda} \times \mathbf D(0, \eta^+)) \ar@{^{(}->}[u]
\\
]]U[[ \ar[r]^-\simeq \ar@{^{(}->}[u]& ]U[ \times \mathbf D(0, 1^-) \ar@{^{(}->}[u].
}
$$

Of course, if $T \subset X_{k}$ is a finite closed subset, we may always assume that $T \subset U$, and we obtain $]]T[[ \simeq ]T[ \times \mathbf D(0, 1^-)$.
If $x$ is a closed point with unramified lifting $a$, we even get
$$
]]x[[ \simeq \mathbf D_{K(a)}(0, 1^-) \times \mathbf D(0, 1^-) \simeq \mathbf B^2_{K(a)}(0, 1^-)
$$
which can also be proved directly as in proposition \ref {pointdisc}.
Finally, we can describe the inverse image by specialization $]]\xi[[$ of the generic point $\xi$ of $X_{k}$: we have
$$
]]\xi[[ \simeq \mathbf D_{]\xi[}(0, 1^-)
$$
where the disc is over the completion $\mathcal H(]\xi[)$ of the field $\mathcal O_{X_{K}^{\mathrm{an}}, ]\xi[}$.

Following Berthelot, we will define now the overconvergence condition.
Note before that the natural inclusion of the first infinitesimal neighborhood $P \hookrightarrow X \times X$, induces by functoriality a morphism $P_{K}^{\mathrm{an}} \hookrightarrow ]]X_{k}[[$ and the projections $q_{1}, q_{2} :P_{K}^{\mathrm{an}} \to X_{K}^{\mathrm{an}}$ introduced in the previous section are induced by $p_{1}$ and $p_{2}$.

\begin{dfn}
Let $Y$ be a locally closed subset of $X_{k}$.
A connection on an $i_{Y}^{-1}\mathcal O_{X_{K}^{\mathrm{an}}}$-module $E$ is \emph{overconvergent} if the corresponding 1-stratification is induced by a \emph{Taylor isomorphism}
$$
\epsilon : i_{Y}^{-1}p_{2}^*i_{Y*}E \simeq i_{Y}^{-1}p_{1}^*i_{Y*} E
$$
on $]]Y[[$ which satisfies the \emph{cocycle condition} (see below) on triple products.
If $Y$ is a closed subset, the connection is just said to be \emph{convergent}.
We will also say that the $\nabla$-module $E$ is \emph{(over-) convergent}.
\end{dfn}

Again, this definition applies in very general geometric situations.
It is also important to insist on the fact that there are no finiteness conditions at that point.
Note that the convergence condition is a lot simpler to describe: if $Z$ is a closed subset in $X_{k}$ (either finite or $X_{k}$ itself), and if we still denote by $p_{1}, p_{2} : ]]Z[[ \to ]Z[$, the maps induced by the projections, we simply have $i_{Y}^{-1}p_{j}^*i_{Y*}E = p_{j}^{*}E$ for $i=1,2$ and the Taylor isomorphism reads
$$
\epsilon : p_{2}^* E \simeq p_{1}^* E
$$
on $]]Z[[$.

It is quite easy to make the \emph{cocycle condition} precise.
We embed diagonally $X$ into the triple product $X \times X \times X$.
We denote by $]]]Y[[[$ the tube of $Y \subset X_{k}$ inside $X_{K}^{\mathrm{an}} \times X_{K}^{\mathrm{an}} \times X_{K}^{\mathrm{an}}$ and by $p_{12}, p_{23}, p_{13} : ]]]X_{k}[[[ \to ]]X_{k}[[$ the maps induced by the projections.
Then the cocycle condition reads
$$
(i_{Y}^{-1}p_{12}^*i_{Y*})(\epsilon) \circ (i_{Y}^{-1}p_{23}^*i_{Y*})(\epsilon) = (i_{Y}^{-1}p_{13}^*i_{Y*})(\epsilon).
$$

\begin{lem}\label{pulbac}
Overconvergence is preserved under restriction to a locally closed subset $Y' \subset Y$.
\end{lem}

\textbf {Proof: }
We denote by the same letter $i$ both inclusion maps $]Y'[ \hookrightarrow ]Y[$ and $]]Y'[[ \hookrightarrow ]]Y[[$.
Then we have for $j=1, 2$,
$$
i_{Y'}^{-1}p_{j}^*i_{Y'*}i^{-1}E = i_{Y'}^{-1}(\mathcal O_{]]X_{k}[[} \otimes_{p_{j}^{-1}\mathcal O_{X_{K}^{\mathrm{an}}}}p_{j}^{-1}i_{Y'*}i^{-1}E)
$$
$$
= i_{Y'}^{-1}\mathcal O_{]]X_{k}[[} \otimes_{i_{Y'}^{-1}p_{j}^{-1}\mathcal O_{X_{K}^{\mathrm{an}}}} i_{Y'}^{-1}p_{j}^{-1}i_{Y'*}i^{-1}E =  i_{Y'}^{-1}\mathcal O_{]]X_{k}[[} \otimes_{i_{Y'}^{-1}p_{j}^{-1}\mathcal O_{X_{K}^{\mathrm{an}}}} i^{-1}p_{j}^{-1}i_{Y*}E
$$
$$
= i^{-1}(i_{Y}^{-1}\mathcal O_{]]X_{k}[[} \otimes_{i_{Y}^{-1}p_{j}^{-1}\mathcal O_{X_{K}^{\mathrm{an}}}} p_{j}^{-1}i_{Y*}E) = i^{-1}(i_{Y}^{-1}p_{j}^*i_{Y*}E).
\quad \Box$$

We will use this explicit isomorphism below for the restriction from $X_{k}$ to some $Y$ in which case it says that the adjunction map is an isomorphism
$$
i_{Y}^{-1}p_{j}^*E \simeq i_{Y}^{-1}p_{j}^*i_{Y*}i_{Y}^{-1}E.
$$

\begin{prop}
If $U$ is an open subset of $X_{k}$ and $E$ is an overconvergent $\nabla$-module on $]U[$, then $i_{U*}E$ is convergent.
If $Z$ is a closed subset of $X_{k}$ and $E$ is a convergent $\nabla$-module on $]Z[$, then $i_{Z!}E$ is convergent.
\end{prop}

\textbf {Proof: }
We want to show that for $j = 1, 2$, we have
$$
p_{j}^*i_{U*} E \simeq i_{U*}i_{U}^{-1}p_{j}^*i_{U*}E.
$$
It is actually sufficient to show that
$$
p_{j}^{-1}i_{U*} E \simeq i_{U*}i_{U}^{-1}p_{j}^{-1}i_{U*}E.
$$
And this follows from the fact that there is a cartesian diagram
$$
\xymatrix{]]U[[ \ar@{^{(}->}[r] \ar[d] & ]]X_{k}[[ \ar[d]^{p_{2}} \\ ]U[ \ar@{^{(}->}[r] & X_{K}^{\mathrm{an}}}
$$
where the horizontal arrows are closed immersions.
Therefore, the 1-stratification extends canonically.

The same type of argument shows that
$$
p_{j}^*i_{Z!} E \simeq i_{Z!}i_{Z}^{-1}p_{j}^*i_{Z!}E
$$
in the second case.
$\quad \Box$

The next proposition will be the first illustration of the power of the overconvergence condition.
Recall that a $\nabla$-module on an analytic variety is said to be \emph{(locally) trivial} if it is (locally) generated by a finite set of horizontal sections.
Note that a locally trivial $\nabla$-module on a disc is always trivial.
For further use, we also state a lemma.

\begin{lem}
If $D \subset X$ is a smooth divisor with reduction $Z$, then the inclusion $D_{K} \hookrightarrow ]Z[$ and its retraction induce an equivalence between coherent $\mathcal O_{D_{K}}$-module and locally trivial $\nabla$-modules on $]Z[$.
\end{lem}

\textbf {Proof: }
We may assume that $Z$ is reduced to one point $x$ with unramified lifting $a$, in which case we are simply considering the inclusion of the origin $\{0\} \hookrightarrow D_{K(a)}(0, 1^-)$ into the disc.
$\quad \Box$

\begin{prop} \label{pointtriv}
If $Z \subset X_{k}$ is a finite closed subset, then any coherent convergent $\nabla$-module on $]Z[$ is locally trivial.
\end{prop}

\textbf {Proof: }
We may assume that $Z$ is reduced to a rational point $x$ and lift it to a rational point $a \in X_{K}$.
Think of $a$ as the ``center'' of $]x[$.
We set $V := a^*E$ which is a finite dimensional vector space.
Next, we consider the composition $X \to \mathrm{Spec}\; \mathcal V \hookrightarrow X$ of the projection $p$ and the section $a$.
Its graph $X \to X \times X$ induces a morphism $\gamma : ]x[ \hookrightarrow ]]x[[$ along which we can pull back the Taylor isomorphism.
We obtain an isomorphism
$$
\phi : \mathcal O_{]x[} \otimes_{K} V = p^* a^* E =  \gamma^*p_{2}^*E \simeq \gamma^*p_{1}^*E = E.
$$
Using the cocycle condition, one sees that the 1-stratification of $E$ is compatible with the trivial one on the left hand side.
More precisely, we use the map $X \times X \to X \times X \times X$ obtained by tensorizing the identity with the above graph.
We may then pull the cocycle condition back along the induced map $]]x[[ \to ]]]x[[[$ and get $\epsilon \circ p_{2}^*(\phi) = p_{1}^*(\phi)$.
$\quad \Box$

Actually, the convergence condition on $\{x\}$ is also called the \emph{Robba condition} and it is a classical result that a finite $\nabla$-module that satisfies the Robba condition on an open disc is automatically trivial.

\begin{cor}
If $E$ is a constructible convergent $\nabla$-module on $X_{K}^{\mathrm{an}}$ and $x \in X_{k}$, then the restriction of $E$ to $]x[$ is trivial. $\quad \Box$
\end{cor}

We now turn to the description of overconvergence on open subsets.
Before doing anything else and although we will continue to prove most statements without referring to rigid cohomology, we should mention the following comparison result:

\begin{prop} \label{isoc}
If $U$ is an open subset of $X_{k}$, the functor $E \mapsto (i_{U*}E)_{0}$ induces an equivalence between the category of coherent overconvergent $\nabla$-modules on $]U[$ and overconvergent isocrystals on $U$.
\end{prop}

\textbf {Proof: }
Using proposition 7.2.13 (and definition 7.2.10) of \cite{LeStum07}, this is simply a translation of the above overconvergence condition into the language of rigid geometry as in proposition \ref{comprig}.
$\quad \Box$

Assume that $t$ is a local parameter defined on some affine open subscheme $\mathrm{Spec}\; A \subset X$ with reduction $U$.
Write as usual $V_{\lambda} = X_{K}^{\mathrm{an}} \setminus ]Z[_{\lambda} = \mathcal M(A_{\lambda})$ where $Z$ is the closed complement of  $U$.
If $M$ is a finite $\nabla$-$A^\dagger_{K}$-module, it extends as usual to some finite $\nabla$-$A_{\lambda}$ module $M_{\lambda}$ with $\lambda < 1$ and we may define the \emph{$\lambda$-radius of convergence} of $M$ as
$$
R(M, \lambda) = \inf_{s \in M_{\lambda}} \inf \left\{\quad \lambda \quad, \quad \underline{\lim}_{k} \left \|\frac 1{k!} \frac{\partial^k}{\partial t^k}(s) \right \|_{\lambda}^{-\frac 1k} \right\}
$$
where $\| - \|_{\lambda}$ is a Banach norm on $M_{\lambda}$.
We say that $M$ is \emph{overconvergent} if
$$
\lim_{\lambda \to 1} R(M, \lambda) = 1.
$$
Those who are interested in differential equations should notice that this condition is equivalent to requiring all the Robba fibers $\mathcal R_{x}(M)$ to be solvable whenever $x \in Z$ (use the ``isometry'' $A_{K}^\dagger \hookrightarrow \mathcal R_{Z}$).
If one is willing to use rigid analytic geometry, the next proposition is a particular case of proposition 7.2.15 of \cite{LeStum07}.

\begin{prop} \label{moddag}
Let $U \subset X_{k}$ be an affine open subset of $X_{k}$.
Let $t$ be a local parameter defined on an algebraic lifting $\mathrm{Spec}\; A$ of $U$.
Then, the functor $\Gamma(]U[, -)$ induces an equivalence between coherent overconvergent $\nabla$-modules on $U$ and finite overconvergent $\nabla$-$A^\dagger_{K}$-modules.
\end{prop}

\textbf {Proof: }
Let $E$ be a $\nabla$-module on $]U[$ and $M$ the corresponding $\nabla$-$A^\dagger_{K}$-module.
One can check that the Taylor isomorphism is necessarily induced by the Taylor series
$$
s \mapsto \sum \frac 1{k!} \frac{\partial^k }{\partial t^k}(s) \tau^k
$$
(see the remark after proposition \ref{eqstrat} below).
The overconvergence condition on $]U[$ means that this series converges on some neighborhood of $]]U[[$ inside $]]X_{k}[[$.
And we may assume that this neighborhood has the form
$$
\cup_{\lambda, \eta} (V_{\lambda} \times \mathbf D(0, \eta^+)) \quad \mathrm{with} \quad \lambda, \eta \to 1
$$
with $V_{\lambda}$ as usual.
The overconvergence condition for $M$ is then the direct translation of the convergence of the series on this specific neighborhood.
$\quad \Box$

We will need below general stratifications using all infinitesimal neighborhoods of $X$.
We denote by $P^{(n)}$ the n-th infinitesimal neighborhood of $X$ in $X \times X$ defined by $\mathcal I^{n+1}$ if $X$ is defined by $\mathcal I$ into $X \times X$, and by $p_{1}^{(n)}, p_{2}^{(n)} = P_{K}^{(n),\mathrm{an}} \to X_{K}^{\mathrm{an}}$ the maps induced by the projections.
With our former notations, we have $P = P^{(1)}$ and $q_{j} = p_{j}^{(1)}$.

\begin{dfn}
If $Y \subset X_{k}$ is a locally closed subset, a \emph{stratification} on an $i_{Y}^{-1}\mathcal O_{X_{K}^{\mathrm{an}}}$-module $E$ is a compatible family of isomorphisms
$$
\epsilon_{n} : i_{Y}^{-1}p_{2}^{(n)*}i_{Y*}E \simeq i_{Y}^{-1}p_{1}^{(n)*}i_{Y*} E
$$
such that the cocycle condition holds on triple products and $\epsilon_{0}$ is the identity.
A morphism of such is an $i_{Y}^{-1}\mathcal O_{X_{K}^{\mathrm{an}}}$-linear map that commutes with the stratifications.
\end{dfn}

\begin{prop}\label{eqstrat}
The category of stratified $i_{Y}^{-1}\mathcal O_{X_{K}^{\mathrm{an}}}$-modules is equivalent the category of $\nabla$-modules on $Y$.
\end{prop}

\textbf {Proof: }
Standard.
$\quad \Box$

It should also be mentioned that the stratification is automatically induced by the Taylor isomorphism when $E$ is overconvergent.

\begin{prop} \label{trulem}
Let $E$ be a $\nabla$-module on $X_{K}^{\mathrm{an}}$, $U$ an affine open subset of $X$ and $Z$ its closed complement.
Assume that $i_{U}^{-1}E$ is a coherent overconvergent $\nabla$-module and that $i_{Z}^{-1}E$ is a locally trivial $\nabla$-module.
Then, $E$ is convergent.
\end{prop}

\textbf {Proof: }
The point is to show that the stratification is induced by a Taylor isomorphism $\epsilon : p_{2}^*E \simeq p_{1}^*E$ on $]]X_{k}[[$.
We already have isomorphisms
$$
\epsilon_{U} : i_{U}^{-1}p_{2}^*E \simeq i_{U}^{-1}p_{1}^*E \quad \mathrm{and} \quad \epsilon_{Z} : i_{Z}^{-1}p_{2}^*E \simeq i_{Z}^{-1}p_{1}^*E
$$
coming from the overconvergence of $E$ on $U$ and $Z$ (see the remark following lemma \ref{pulbac}).
And we can use again the fact that a sheaf on a topological space is uniquely determined by its restriction to an open subset, its restriction to the closed complement and the adjunction map.
It is therefore sufficient to show that the diagram
$$
\xymatrix{
i_{U}^{-1}p_{2}^*E \ar[r] \ar[d]^{\epsilon_{U}} & i_{U}^{-1}i_{Z*}i_{Z}^{-1}p_{2}^*E \ar[d]^{\epsilon_{Z}}
\\
i_{U}^{-1}p_{1}^*E \ar[r] & i_{U}^{-1}i_{Z*}i_{Z}^{-1}p_{1}^*E}
$$
is commutative.
Since the analogous diagram with stratifications is commutative by hypothesis (we have a stratification on $E$), it is sufficient to prove that the canonical map
$$
\mathrm{Hom}_{i_{U}^{-1}\mathcal O_{]]X_{k}[[}}(i_{U}^{-1}p_{2}^*E , i_{U}^{-1}i_{Z*}i_{Z}^{-1}p_{1}^*E) \to \varprojlim \mathrm{Hom}_{i_{U}^{-1}\mathcal O_{P_{K}^{(n),\mathrm{an}}}}(i_{U}^{-1}p_{2}^{(n)*}E , i_{U}^{-1}i_{Z*}i_{Z}^{-1}p_{1}^{(n)*}E) 
$$
is injective.
Since $i_{U}^{-1}E$ is finitely presented and $i_{Z}^{-1}E$ is trivial, it is sufficient to prove that
the canonical map
$$
\Gamma(]]U[[, i_{U}^{-1}i_{Z*}i_{Z}^{-1}p_{1}^*\mathcal O_{]]X_{k}[[}) \to \varprojlim 
\Gamma(]U[, i_{U}^{-1}i_{Z*}i_{Z}^{-1}p_{1}^{(n)*}  \mathcal O_{P_{K}^{(n),\mathrm{an}}}) 
$$
is injective.
By considering the stalks and since $]]\xi[[$ is closed in $]]X_{k}[[$, one sees that
$$
i_{U}^{-1}i_{Z*}i_{Z}^{-1}p_{1}^*\mathcal O_{]]X_{k}[[} = i_{U}^{-1}i_{Z*}\mathcal O_{]]Z[[} = i_{\xi*}i_{\xi}^{-1}i_{Z*}\mathcal O_{]]Z[[}.
$$
Also, if we fix a local parameter $t$ on $X$, we can easily identify the right hand side with $\mathcal R_{Z}[[\tau]]$.
We are led to check that the canonical map
$$
\Gamma(]]\xi[[, i_{\xi}^{-1}i_{Z*}\mathcal O_{]]Z[[}) \to \mathcal R_{Z}[[\tau]]
$$
is injective.

We may assume that the local parameter $t$ is defined on a an affine open subscheme $\mathrm{Spec}\; A$ whose reduction $U$ contains $Z$.
As usual, if $T$ denotes the closed complement of $U$, we let $V_{\lambda} := X_{K}^{\mathrm{an}} \setminus ]T[_{\lambda}$.
It is sufficient to show that if $V'$ is a neighborhood of $]]\xi[[$ inside $]]X_{k}[[$, the map
$$
\Gamma(V' \cap ]Z[, \mathcal O_{V'}) \to \mathcal R_{Z}[[\tau]]
$$
is injective.
After removing some points in $U$ if necessary, we may assume that
$$
V' = \cup_{\lambda, \eta} (V_{\lambda} \times \mathbf D(0, \eta^+))
$$
with $\lambda, \eta \to 1$.
We are then reduced to showing that the map
$$
\Gamma(]Z[ \setminus ]Z[_{\lambda} \times \mathbf D(0, \eta^+), \mathcal O) \to \mathcal R_{Z}[[\tau]]
$$
is injective.
An it is sufficient to consider the obvious injective map
$$
\Gamma([Z]_{\mu} \setminus ]Z[_{\lambda} \times \mathbf D(0, \eta^+), \mathcal O) \hookrightarrow \Gamma([Z]_{\mu} \setminus ]Z[_{\lambda}, \mathcal O)[[\tau]]
$$
and take inverse limit when $\mu \to 1$.
$\quad \Box$

\begin{cor} \label{locover}
Let $E$ be a constructible $\nabla$-module on $X_{K}^{\mathrm{an}}$.
The connection on $E$ is convergent if and only if there exists a finite covering of $X_{k}$ by locally closed subsets $Y$ such that the connection is overconvergent on each $Y$.
\end{cor}

\textbf {Proof: }
Follows also from proposition \ref{trulem}.
$\quad \Box$

\begin{cor} \label{unscrew}
A $\nabla$-module $E$ on $X_{K}^{\mathrm{an}}$ is constructible convergent if and only if there exists an exact sequence
$$
0 \to i_{Z!}E_{Z} \to E \to i_{U*} E_{U} \to 0 
$$
where $U$ is a non-empty affine open subset of $X_{k}$ with closed complement $Z$, $E_{U}$ is a coherent overconvergent $\nabla$-module on $]U[$ and $E_{Z}$ is a locally trivial $\nabla$-module on $]Z[$.
We can even assume that $E_{U}$ is free.
\end{cor}

\textbf {Proof: }
Follows from propositions \ref{exactseq} and \ref{trulem}.
$\quad \Box$

Recall that the notion of $\nabla$-extension was introduced in definition \ref{nablext}.

\begin{dfn}
A $\nabla$-extension $\mathcal R(M) \to \mathcal R(\mathcal F)$ is said to be \emph{convergent} if $M$ is overconvergent and $\mathcal F$ is locally trivial.
\end{dfn}

In other words, a convergent $\nabla$-extension is given by a an overconvergent $\nabla$-$A^\dagger_{K}$-module $M$ where $\mathrm{Spec}\; A \subset X$ is non empty affine open subset with reduction $U$, a collection of finite dimensional $K(a)$-vector spaces $H_{a}$ for each point $x$ not in $U$, where $a$ is an unramified lifting of $x$, and horizontal $A^\dagger_{K}$-linear maps $M \to \mathcal R_{a} \otimes_{K(a)} H_{a}$.

\begin{thm} \label{convext}
The category of constructible convergent $\nabla$-modules on $X_{K}^{\mathrm{an}}$ is equivalent to the category of convergent $\nabla$-extensions modulo equivalence.
\end{thm}

\textbf {Proof: }
Follows from propositions \ref{equinab} and \ref{trulem}.
$\quad \Box$

Finally, we study the relation between overconvergence and Frobenius:

\begin{prop} 
Assume that $F: X \to X$ is a lifting of the absolute Frobenius of $X_{k}$.
 If $Y$ is a locally closed subset of $X_{k}$ and $E$ is an overconvergent $\nabla$-module on $]Y[$, then $i_{Y}^{-1}F_{K}^{\mathrm{an}*}i_{Y*}E$ is also overconvergent.
\end{prop}

\textbf {Proof: }
If we denote by $F_{]]X_{k}[[}$ the map induced on $]]X_{k}[[$ by $F \times F$, then for $j = 1, 2$, we have
$$
i_{Y}^{-1}p_{j}^*i_{Y*}i_{Y}^{-1} F_{K}^{\mathrm{an}*} i_{Y*} = i_{Y}^{-1}F_{]]X_{k}[[}^*i_{Y*}i_{Y}^{-1}p_{j}^*i_{Y*}
$$
and the result follows.
$\quad \Box$

\begin{prop}
If $F: X \to X$ is a lifting of the absolute Frobenius of $X_{k}$, then any constructible $F$-$\nabla$-module on $X_{K}^{\mathrm{an}}$ is convergent.
\end{prop}

\textbf {Proof: }
This is Dwork's trick.
Using proposition \ref{trulem}, we may assume either that $E = i_{U*}E_{U}$ with $U \subset X_{k}$ affine open and $E_{U}$ a coherent $i_{U}^{-1}\mathcal O_{X_{K}^{\mathrm{an}}}$-module or $E = i_{x!}E_{x}$ with $x \in X_{k}$ closed and $E_{x}$ coherent on $]x[$.
The first case follows from proposition \ref{moddag} using for example corollary 8.3.9 of \cite{LeStum07}.
The second case is the classical fact that an (finite) $F$-$\nabla$-module on an open disc is always trivial.
$\quad \Box$

\section{Constructibility and $\mathcal D$-modules}

We introduce now the ring of arithmetic differential operators on $\hat X$.
It can be described as follows.
First of all, if $\mathcal D_{\hat X}$ denotes the sheaf of algebraic differential operators on $\hat X$ as before, and $\widehat{\mathcal D}_{\hat X}$ is its $p$-adic completion, we have
$$
\mathcal D_{\hat X} \subset \mathcal D^\dagger_{\hat X} \subset \widehat{\mathcal D}_{\hat X}.
$$
Moreover, as shown in proposition 2.4.4 of \cite{Berthelot96}, if  $t$ is a local parameter on $X$ defined on some subset $\mathcal U = \mathrm{Spf}\; \mathcal A \subset \hat X$, then
$$
\Gamma(\mathcal U, \mathcal D^\dagger_{\hat X}) = \{\sum_{k=0}^{\infty} \frac 1{k!} f_{k} \frac {\partial^k}{\partial t}, \quad \exists c > 0, \eta < 1, \quad \|f_{k}\| \leq c\eta^k\}.
$$

We are interested in $\mathcal D^\dagger_{\hat X\mathbf Q}$-modules.
The case of coherent $\mathcal D^\dagger_{\hat X \mathbf Q}$-modules with finite support is easy to settle.
If $D \subset X$ is a smooth divisor with reduction $Z$, then the direct image functor
$$
\xymatrix@R0cm{
D^b_{\mathrm{coh}}(\mathcal O_{\hat D\mathbf Q}) \ar[rr]^{i_{Z+}} && D^b_{\mathrm{coh}}(\mathcal D^\dagger_{\hat X\mathbf Q})
\\
H \ar@{|->}[rr] && i_{Z*}(i^{*}_{Z} \mathcal Hom_{\mathcal O_{\hat X}}(\Omega^1_{\hat X}, \mathcal D^\dagger_{\hat X})_{\mathbf Q} \otimes_{\mathcal O_{\hat D \mathbf Q}} H)
}
$$
induces an equivalence between coherent $\mathcal O_{\hat D\mathbf Q}$-modules and coherent $\mathcal D^\dagger_{\hat X\mathbf Q}$-modules with support in $Z$ (see section 5.3.3 of \cite{Berthelot02} for the general statement).
Its inverse is induced by the exceptional inverse image $i_{Z}^{!}$.

Recall that if $Z$ is a finite closed subset of $X$, there is a short exact sequence
$$
0 \to \mathcal O_{\widehat{X}} \to \mathcal O_{\widehat{X}}({}^\dagger Z) \to \mathcal H^\dagger_Z  \to 0
$$
of $\mathcal D^\dagger_{\hat X}$-module.
Moreover, $\mathcal O_{\hat{X}\mathbf Q}$, $\mathcal O_{\widehat{X}}({}^\dagger Z)_{\mathbf Q}$ and $\mathcal H^\dagger_{Z\mathbf Q} $ are all $\mathcal D^\dagger_{\hat X\mathbf Q}$-coherent.

The next result is also well-known (formula 4.4.5.2 of \cite{Berthelot02}) but rather easy to check (and instructive) in our situation.

\begin{prop}\label{compar}
If $D \subset X$ is a smooth divisor with reduction $Z$, and $H$ is a coherent $\mathcal O_{\hat D \mathbf Q}$-module, there is a canonical isomorphism of $\mathcal D^\dagger_{\hat X\mathbf Q}$-modules
$$
i_{Z+} H \simeq \mathcal H^\dagger_{Z\mathbf Q} \otimes_{\mathcal O_{\hat D \mathbf Q}} H.
$$
\end{prop}

\textbf {Proof: }
We may assume that $Z = \{x\}$ where $x$ is the only zero of a local parameter $t$ defined on some formal affine open subset $\mathcal U = \mathrm{Spf} \mathcal A$, and that $H = K(a)$ with $a$ an unramified lifting of $x$.
We have a commutative diagram with exact rows (for multiplication on the left)
$$
\xymatrix{
0 \ar[r] & \Gamma(\mathcal U, \mathcal D^\dagger_{\hat X}) \ar[rr]^{\frac \partial {\partial t}} \ar@{=}[d] && \Gamma(\mathcal U, \mathcal D^\dagger_{\hat X})  \ar[rr] \ar[d]^t && \mathcal A \ar[r] \ar[d] & 0
\\
0 \ar[r] & \Gamma(\mathcal U, \mathcal D^\dagger_{\hat X}) \ar[rr]^{t\frac \partial {\partial t} + 1} && \Gamma(\mathcal U, \mathcal D^\dagger_{\hat X})  \ar[rr]^{1 \mapsto \frac 1t} && \mathcal A[1/t]^{\dagger} \ar[r] & 0
}
$$
from which we deduce an isomorphism
$$
\Gamma(\mathcal U, \mathcal D^\dagger_{\hat X})/\Gamma(\mathcal U, \mathcal D^\dagger_{\hat X})t \simeq \mathcal A[1/t]^\dagger/\mathcal A.
$$
The right hand side is exactly $\Gamma(\mathcal U,\mathcal H^\dagger_Z)$.
On the other hand, we have
$$
\Gamma(\mathcal U, i_{x+}K(a)) = K(a) \otimes_{\mathcal A_{K}}\mathrm{Hom}_{\mathcal A_{K}}(\mathcal A_{K}\mathrm dt, \Gamma(\mathcal U, \mathcal D^\dagger_{\hat X\mathbf Q}) \simeq \Gamma(\mathcal U, \mathcal D^\dagger_{\hat X\mathbf Q})/\Gamma(\mathcal U, \mathcal D^\dagger_{\hat X \mathbf Q})t. 
$$
We get an isomorphism as expected which is easily seen to be independent of the choices (just multiply $t$ by an invertible element in $\mathcal A_{K}$).
$\quad \Box$

The next proposition is a very fancy way of stating an almost trivial result.
However, it may be seen as an analog of theorem 6.5 of \cite{Kashiwara04} in a very simple situation.
This may also be seen as a special case of theorem 2.5.10 of \cite{Caro09b}.

\begin{prop} \label{fineq}
If $Z$ is a finite closed subset of $X_{k}$, then $\mathrm R \widetilde{\mathrm{sp}}_{*} i_{Z!}$ induces an equivalence between the category of constructible convergent $\nabla$-modules on $]Z[$ and the category complexes of $\mathcal D^\dagger_{\hat X\mathbf Q}$-modules with support in $Z$ and coherent cohomology concentrated in degree $1$.
\end{prop}

\textbf {Proof: }
If $D \subset X$ is a smooth divisor with reduction $Z$, and $H$ is a coherent $\mathcal O_{\hat D \mathbf Q}$-module, it follows from proposition \ref{compar} and proposition \ref{pointwis} that
$$
\mathrm R \widetilde{\mathrm{sp}}_{*} i_{Z!}(\mathcal O_{]Z[} \otimes_{\mathcal O_{\hat D \mathbf Q}} H) = i_{Z+}H[1].
$$
Moreover, we know that $\mathcal O_{]Z[} \otimes_{\mathcal O_{\hat D}} - $ is an equivalence of categories between coherent $\mathcal O_{\hat D_{K}}$-modules and locally trivial $\nabla$-modules on $]Z[$, and that $i_{Z+}$ is an equivalence of categories between coherent $\mathcal O_{\hat D\mathbf Q}$-modules and coherent $\mathcal D^\dagger_{\hat X\mathbf Q}$-modules with support in $Z$.
Since specialization $\hat D_{K} \to Z$ is bijective, we can identify $\mathcal O_{\hat D\mathbf Q}$ with $\mathcal O_{\hat D_{K}}$ and consequently $\mathcal O_{\hat D\mathbf Q}$-modules with $\mathcal O_{\hat D_{K}}$-modules.
Finally, we know from proposition \ref{pointtriv} that constructible convergent $\nabla$-modules on $]Z[$ are locally trivial (and conversely).
$\quad \Box$

We now turn to the case of an affine open subset of $X_{k}$.
It is necessary to introduce the ring of arithmetic differential operators $\mathcal D^\dagger_{\hat X}({}^\dagger Z)$  with overconvergent poles along a non empty finite closed subset $Z$.
Since we are only interested in this ring modulo torsion, it is more convenient to use the modified definition given in 2.6.2 of \cite{Huyghe03}.

Just before proposition \ref{weakequi}, we introduced the ring $\mathcal D_{\hat X}({}^\dagger Z)$ of algebraic differential operators with overconvergent poles.
If $j : X_{k} \setminus Z \hookrightarrow X_{k}$ denotes the inclusion map, we will have
$$
\mathrm{both}\ \mathcal D^\dagger_{\hat X} \ \mathrm{and} \ \mathcal D_{\hat X}({}^\dagger Z) \subset \mathcal D^\dagger_{\hat X}({}^\dagger Z) \subset j_{*}j^{-1}\mathcal D^\dagger_{\hat X}.
$$
If there exists a local parameter $t$ defined on some open subset $\mathcal U = \mathrm{Spf}\; \mathcal A$ such that $Z \cap \mathcal U = \{x\}$ where $x$ is the only zero of $t$, then
$$
\Gamma(\mathcal U, \mathcal D^\dagger_{\hat X}({}^\dagger Z)) = \{\sum_{k,l=0}^{\infty} \frac 1{k!} \frac {f_{k,l}}{t^l} \frac {\partial^k}{\partial t^k}, \quad \exists c > 0, \eta < 1, \quad \|f_{k,l}\| \leq c\eta^{k+l}\}.
$$

Note that $\mathcal D^\dagger_{\hat X}({}^\dagger Z)_{\mathbf Q}$-modules that are $\mathcal O_{\hat X\mathbf Q}({}^\dagger Z)$-coherent form a full subcategory of the category of $\mathcal D_{\hat X}({}^\dagger Z)_{\mathbf Q}$ modules.
In other words, the forgetful functor is fully faithful (use Spencer complexes to show this).

We will have to consider scalar extension and we will write when $\mathcal E$ is a coherent $\mathcal D^\dagger_{\hat X\mathbf Q}$-module,
$$
\mathcal E({}^\dagger Z) := \mathcal D^\dagger_{\hat X}({}^\dagger Z)_{\mathbf Q} \otimes_{\mathcal D^\dagger_{\hat X\mathbf Q}} \mathcal E.
$$
Then, there exists an exact triangle
$$
i_{Z+}i_{Z}^{!} \mathcal E  \to \mathcal E \to \mathcal E({}^\dagger Z) \to .
$$

We want to explain now, as in proposition 4.4.3 of \cite{Berthelot96}, that if $U = X \setminus Z$, and $E$ is a coherent overconvergent $\nabla$-module on $]U[$, then the action of $\mathcal D_{\hat X}({}^\dagger Z)$ on $\widetilde{\mathrm{sp}}_{*}i_{U*}E$ extends to an action of  $\mathcal D^\dagger_{\hat X}({}^\dagger Z)$.

Assume that $t$ is a local parameter defined on some open subset $\mathcal U' = \mathrm{Spf}\; \mathcal A'$ with $Z \cap \mathcal U' = \{x\}$ where $x$ is the only zero of $t$ and that $t$ is defined on a lifting $\mathrm{Spec}\, A$ of $U$.
The coherent $\nabla$-module $E$ is given by a finite $\nabla$-$A_{K}^\dagger$-module $M$ and
$$
\Gamma(\mathcal U', \widetilde{\mathrm{sp}}_*i_{U*} E) = \mathcal A'[1/t]^\dagger_{\mathbf Q} \otimes_{A^\dagger_{K}} M =: M'.
$$
Also, we can write $M = \varinjlim_{\lambda < 1} M_{\lambda}$ where $M_{\lambda}$ is a finite $\nabla$-$A_{\lambda}$-module so that
$$
M' = \varinjlim_{\lambda < 1} M'_{\lambda} \quad \mathrm{with} \quad M'_{\lambda} := \mathcal A'_{K}\{\lambda/t\} \otimes_{A_{\lambda}} M.
$$
If $P \in \Gamma(\mathcal U', \mathcal D^\dagger_{\hat X}({}^\dagger Z))$, then there exists $c > 0$ and $\eta < 1$ such that
$$
P = \sum_{k,l=0}^{\infty}  \frac 1{k!} \frac {f_{k,l}}{t^l} \frac {\partial^k}{\partial t^k} \quad \mathrm{with} \quad \|f_{k,l}\| \leq c\eta^{k+l}.
$$
If $M$ is overconvergent and $\lambda$ close enough to $1$, we have $\|\frac 1{k!}\frac {\partial^k(s)}{\partial t^k}\|\eta^k \to 0$ for $s \in M_{\lambda}$.
Of course, we will also have $\lambda > \eta$ for $\lambda$ close enough to one.
Then, we see that if $s'$ denotes the restriction of $s$ to $\mathcal U'_{K} \cap V_{\lambda}$, we have
$$
\left\| \frac 1{k!} \frac {f_{k,l}}{t^l} \frac {\partial^k(s')}{\partial t^k}\right\| \leq c \frac {\eta^{k+l} }{\lambda^l} \left\|\frac 1{k!} \frac {\partial^k(s')}{\partial t^k} \right \| \leq c \left \|\frac 1{k!} \frac {\partial^k(s')}{\partial t^k} \right\|\eta^k \to 0
$$
and we can set
$$
P(s') = \sum_{k,l=0}^{\infty}  \frac 1{k!} \frac {f_{k,l}}{t^l} \frac {\partial^k(s')}{\partial t^k}.
$$
In other words, we can extend the action by continuity on $s'$.
Since the question is local on $\hat X$ (and also, that it is sufficient to define the action on a finite set of generators), we see that the action of $\mathcal D_{\hat X}({}^\dagger Z)$ on $\widetilde{\mathrm{sp}}_{*}i_{U*}E$ does extend to an action of $\mathcal D^\dagger_{\hat X}({}^\dagger Z)$.

Of course, by continuity, any horizontal map will give a $\mathcal D^\dagger_{\hat X}({}^\dagger Z)$-linear map and we have proven the following:

\begin{prop} \label{unpub}
If $U$ is a non empty affine open subset of $X_{k}$ and $Z$ denotes its closed complement, then the functor $\widetilde{\mathrm{sp}}_{*}i_{U*}$ induces a fully faithful functor form the category of coherent overconvergent $\nabla$-modules $E$ on $]U[$ to the category of $\mathcal D^\dagger_{\hat X}({}^\dagger Z)_{\mathbf Q}$-modules.$\quad \Box$
\end{prop}

Actually, it is a lot better:

\begin{thm}[Berthelot] \label{unpub}
The essential image of this functor is the category of coherent $\mathcal D^\dagger_{\hat X}({}^\dagger Z)_{\mathbf Q}$-modules that are $\mathcal O_{\hat X}({}^\dagger Z)_{\mathbf Q}$-coherent.
\end{thm}

\textbf {Proof: }
This is exactly the content of a letter from Berthelot to Caro (\cite{Berthelot07}).
Actually, the condition is even weaker since we can simply assume that the restriction to the formal lifting $\mathcal U$ of $U$ is $\mathcal O_{\mathcal U}$-coherent.
$\quad \Box$

We want to extend these results to constructible convergent $\nabla$-modules on $X_{K}^{\mathrm{an}}$.
For the moment, we can prove the following:

\begin{prop}
$\mathrm R\widetilde{\mathrm{sp}}_*$ induces a functor from the category of  constructible convergent $\nabla$-modules on $X_{K}^{\mathrm{an}}$ to the category of perverse $\mathcal D^\dagger_{\hat X\mathbf Q}$-modules.
\end{prop}

\textbf {Proof: }
We know from proposition \ref{nablam} that $\mathrm R\widetilde{\mathrm{sp}}_* E$ is a perverse $\mathcal D_{\hat X\mathbf Q}$-module.
Moreover, we gave in proposition \ref{descdir} a description of the terms of this complex.
Since $\mathcal H^\dagger_Z$ is a $\mathcal D^\dagger_{\hat X}$-module, it follows from proposition \ref{unpub} that both terms of $\mathrm R\widetilde{\mathrm{sp}}_* E$ are $\mathcal D^\dagger_{\hat X\mathbf Q}$-modules.
It only remains to check that the non-trivial map in this complex is $\mathcal D^\dagger_{\hat X\mathbf Q}$-linear.

We showed in theorem \ref{convext} that a constructible convergent $\nabla$-module $E$ on $X_{K}^{\mathrm{an}}$ is given by an overconvergent $\nabla$-$A^\dagger_{K}$-module $M$ where $\mathrm{Spec}\; A \subset X$ is non empty affine open subset with reduction $U$, a collection of finite dimensional $K(a)$-vector spaces $H_{a}$ for each point $x \not \in U$, where $a$ is an unramified lifting of $x$, and horizontal $A^\dagger_{K}$-linear maps $M \to \mathcal R_{a} \otimes_{K(a)} H_{a}$.
We explained at the end of section \ref{secsp} how to construct $\mathrm R\widetilde{\mathrm{sp}}_* E$ from these data.
If $t$ is a local parameter defined on some open subset $\mathcal U' = \mathrm{Spf}\; \mathcal A'$ with $Z \cap \mathcal U' = \{x\}$ where $x$ is the only zero of $t$, and $t$ is actually defined on a lifting $\mathrm{Spec}\, A$ of $U$, we have
$$
\Gamma(U', \mathrm R\widetilde{\mathrm{sp}}_* E) = \left[\mathcal A'[1/t]^\dagger_{\mathbf Q} \otimes_{A^\dagger_{K}} M \to \delta_{a} \otimes_{K(a)} H_{a}\right].
$$
Thus, we must show that the morphisms
$$
\mathcal A'[1/t]^\dagger_{\mathbf Q} \otimes_{A^\dagger_{K}} M \to \delta_{a} \otimes_{K(a)} H_{a}
$$
are $\Gamma(\mathcal U', \mathcal D^\dagger_{\hat X})$-linear.
This follows from the fact that they are horizontal and continuous (recall that the action if defined on both sides by continuity).
$\quad \Box$

\begin{lem} \label{holon}
Assume that $F: X \to X$ is a lifting of the absolute Frobenius of $X_{k}$.
If $E$ is a constructible convergent $\nabla$-module on $X_{K}^{\mathrm{an}}$, then the isomorphism
$$
\hat F^*\mathrm R\widetilde{\mathrm{sp}}_* E \simeq \mathrm R\widetilde{\mathrm{sp}}_* F_{K}^{\mathrm{an}*}E
$$
is $\mathcal D^\dagger_{\hat X\mathbf Q}$-linear.
\end{lem}

\textbf {Proof: }
Again, this follows from the fact that it is horizontal and a continuity argument.
$\quad \Box$

The following is a consequence of the local monodromy theorem:

\begin{prop} \label{holon}
Let $F: X \to X$ be a lifting of the absolute Frobenius of $X_{k}$.
If $E$ is a constructible $F$-$\nabla$-module on $X_{K}^{\mathrm{an}}$, then $\mathrm R\widetilde{\mathrm{sp}}_* E$ is a holonomic $F$-$\mathcal D^\dagger_{\hat X\mathbf Q}$-module.
\end{prop}

\textbf {Proof: }
Thanks to propositions \ref{unscrew} and \ref{fineq}, we may assume that $E = i_{U*}E_{U}$ with $E_{U}$ coherent.
We showed in proposition \ref{sperig} that $\widetilde{\mathrm{sp}}_* i_{U*}E = \mathrm{sp}_{*}(i_{U*}E_{U})_{0}$ and in proposition \ref{isoc} that $(i_{U*}E_{U})_{0}$ is an overconvergent $F$-isocrystal.
Then, it is shown in theorem 4.3.4 of \cite{Caro06} (see also \cite{HuygheTrihan07}) that $\mathrm{sp}_{*}(i_{U*}E_{U})_{0}$ is holonomic. 
$\quad \Box$

However, it is not clear that $\mathrm R\widetilde{\mathrm{sp}}_*$ is fully faithful on constructible $F$-$\nabla$-modules or that its image is  exactly the category of perverse holonomic $F$-$\mathcal D^\dagger_{\hat X\mathbf Q}$-modules.
Some more work is necessary.

\section{Formal fibers}

If $x \subset X_{k}$ is a closed point, we can consider the completion $\hat x$ of $X$ along $x$ (usually written $\widehat X^x$).
Note that if $a$ is an unramified point over $x$, then $\hat x \simeq \mathrm{Spf}\; \mathcal V(a)[[t]]$.
As explained in \cite{Crew06} and \cite{Crew10*}, although $\hat x$ is not a $p$-adic formal scheme, there exists a beautiful theory of arithmetic $\mathcal D$-modules on $\hat x$ which is completely analogous to the theory for $\hat X$.
The ring $\mathcal O_{\hat x \mathbf Q}$ is isomorphic to the set of bounded function on the open unit disc over $K(a)$.
The ring $\mathcal R^{\mathrm{bd}}_{x} := \mathcal O_{\hat x}({}^\dagger x)_{\mathbf Q}$ is isomorphic to the bounded Robba ring $\mathcal R_{a}^{\mathrm{bd}}$ over $K(a)$.
And we have
$$
\delta_{x} = \mathcal R_{x}/\mathcal O_{x}^{\mathrm{an}} \simeq \mathcal R^{\mathrm{bd}}_{x}/\mathcal O_{\hat x\mathbf Q}.
$$
The rings $\mathcal D^\dagger_{\hat x}$ and $\mathcal D^\dagger_{\hat x}({}^\dagger x)$ have exactly  the same description as above (using the Gauss norm).

There is a canonical morphism of formal schemes $i_{\hat x} : \hat x \to \hat X$ which is formally \'etale.
One can define for a coherent $\mathcal D^\dagger_{\hat X\mathbf Q}$-module $\mathcal E$, its \emph{exceptional inverse image}
$$
i_{\hat x}^{!}\mathcal E  := \mathcal D^\dagger_{\hat x\mathbf Q} \otimes_{\mathcal D^\dagger_{\hat X\mathbf Q}} \mathcal E
$$
which is a coherent $\mathcal D^\dagger_{\hat x\mathbf Q}$-module, and we have the following adjunction formula:

\begin{lem} \label{crewtech}
If $\mathcal E$ is a coherent $\mathcal D^\dagger_{\hat X\mathbf Q}$-module (or perfect complex) and $x \in X_{k}$ is a closed point, then
$$
\mathrm R\mathcal Hom_{\mathcal D^\dagger_{\hat X\mathbf Q}}(\mathcal E, \mathcal H^\dagger_{x\mathbf Q}) \simeq i_{x*}\mathrm {RHom}_{\mathcal D^\dagger_{\hat x\mathbf Q}}(i_{\hat x}^{!}\mathcal E, \delta_{x}).
$$
\end{lem}

\textbf {Proof: }
We may clearly assume that $\mathcal E =  \mathcal D^\dagger_{\hat X\mathbf Q}$ and we fall back on lemma \ref{spclo}.
$\quad \Box$

\begin{dfn}
Let $\mathcal E$ be a coherent $\mathcal D^\dagger_{\hat X\mathbf Q}$-module (or perfect complex) and $x \in X_{k}$, a closed point.
Then, the \emph{bounded Robba fiber} of $\mathcal E$ at $x$ is
$$
\mathcal R^{\mathrm{bd}}_{x}(\mathcal E) := i_{\hat x}^{!}\mathcal E({}^\dagger x).
$$
\end{dfn}

If $Z \subset X_{k}$ is a finite closed subset and $x \in Z$, we can also define for a coherent $\mathcal D^\dagger_{\hat X}({}^\dagger Z)_{\mathbf Q}$-module $\mathcal E$, its \emph{exceptional inverse image}
$$
i_{\hat x}^{!}\mathcal E := \mathcal D^\dagger_{\hat x}({}^\dagger x)_{\mathbf Q} \otimes_{\mathcal D^\dagger_{\hat X}({}^\dagger Z)_{\mathbf Q}} \mathcal E.
$$
Note that if $\mathcal E$ is a coherent $\mathcal D^\dagger_{\hat X\mathbf Q}$-module, we have
$$
\mathcal R^{\mathrm{bd}}_{x}(\mathcal E) \simeq i_{\hat x}^{!}\mathcal E({}^\dagger Z).
$$
In particular, the bounded Robba fiber is a generic invariant in the sense that
$$
\mathcal R^{\mathrm{bd}}_{x}(\mathcal E) \simeq \mathcal R^{\mathrm{bd}}_{x}(\mathcal E({}^\dagger Z))
$$
whenever $Z$ is a finite closed subset of $X_{k}$.
Finally, if $\mathcal E$ is a $\mathcal D^\dagger_{\hat X}({}^\dagger Z)_{\mathbf Q}$-module which is coherent both as $\mathcal D^\dagger_{\hat X \mathbf Q}$-module and as $\mathcal O_{\hat X}({}^\dagger Z)_{\mathbf Q}$-module, then
$$
\mathcal R^{\mathrm{bd}}_{x}(\mathcal E) \simeq \mathcal R^{\mathrm{bd}}_{x} \otimes_{\mathcal O_{\hat X}({}^\dagger Z)_{\mathbf Q}} \mathcal E
$$

is a finite $\nabla$-module on $\mathcal R^{\mathrm{bd}}_{x}$.

\begin{dfn} \label{frobatx}
Let $\mathcal E$ be a \emph{coherent} $\mathcal D^\dagger_{\hat X\mathbf Q}$-module.
Then, $E$ has \emph{Frobenius type at a closed point $x$} (resp. has \emph{Frobenius type}) if $\mathcal R^{\mathrm{bd}}_{x}(\mathcal E)$ has a Frobenius structure (resp. a Frobenius structure at all closed points $x \in X_{k}$).
\end{dfn}

Of course, a coherent $F$-$\mathcal D^\dagger_{\hat X\mathbf Q}$-module has Frobenius type.

\begin{prop}
Let $Z \subset X_{k}$ be a finite closed subset and $x \in Z$.
Let $\mathcal E$ be a $\mathcal D^\dagger_{\hat X}({}^\dagger Z)_{\mathbf Q}$-module which is coherent both as $\mathcal D^\dagger_{\hat X \mathbf Q}$-module and as $\mathcal O_{\hat X}({}^\dagger Z)_{\mathbf Q}$-module.
Assume that $\mathcal E$ has Frobenius type at $x$.
Then, we have
$$
\mathrm R\mathcal Hom_{\mathcal D^\dagger_{\hat X\mathbf Q}}(\mathcal E, \mathcal H^\dagger_{x}) \simeq i_{x*} \mathrm {RHom}_{\nabla}(\mathcal R_{x}(\mathcal E), \mathcal R_{x}).
$$
\end{prop}

\textbf {Proof: }
Note first that, in our situation, we have $i_{\hat x}^{!}\mathcal E = \mathcal R^{\mathrm{bd}}_{x}(\mathcal E)$.
Thus, it follows from lemma \ref{crewtech} that
$$
\mathrm R\mathcal Hom_{\mathcal D^\dagger_{\hat X\mathbf Q}}(\mathcal E, \mathcal H^\dagger_{x}) \simeq i_{x*}\mathrm {RHom}_{\mathcal D^\dagger_{\hat x\mathbf Q}}(\mathcal R^{\mathrm{bd}}_{x}(\mathcal E), \delta_{x}).
$$
Since $\mathcal R^{\mathrm{bd}}_{x}(\mathcal E)$ is a finite $F$-$\nabla$-module on $\mathcal R^{\mathrm{bd}}_{x}$, if follows from theorem 3.3 of \cite{Crew06} and proposition 5.6 of \cite{Crew10*} that
$$
\mathrm {RHom}_{\mathcal D^\dagger_{\hat x\mathbf Q}}(\mathcal R^{\mathrm{bd}}_{x}(\mathcal E),  \mathcal O_{x}^{\mathrm{an}}) = 0.
$$
But there is an exact sequence
$$
0 \to \mathcal O_{x}^{\mathrm{an}} \to \mathcal R_{x} \to \delta_{x} \to 0,
$$
and we obtain
$$
\mathrm {RHom}_{\mathcal D^\dagger_{\hat x\mathbf Q}}(\mathcal R^{\mathrm{bd}}_{x}(\mathcal E), \delta_{x}) = \mathrm {RHom}_{\mathcal D^\dagger_{\hat x\mathbf Q}}(\mathcal R^{\mathrm{bd}}_{x}(\mathcal E), \mathcal R_{x})
$$
$$
= \mathrm {RHom}_{\mathcal D^\dagger_{\hat x}({}^\dagger x)_{\mathbf Q}}(\mathcal R^{\mathrm{bd}}_{x}(\mathcal E), \mathcal R_{x}) = \mathrm {RHom}_{\nabla}(\mathcal R^{\mathrm{bd}}_{x}(\mathcal E), \mathcal R_{x})
$$
$$
= \mathrm {RHom}_{\nabla}(\mathcal R_{x}(\mathcal E), \mathcal R_{x})
\quad \Box$$

We consider now the analog definitions on the constructible side.
Recall that if $x \in X_{k}$ is a closed point, we introduced in definition \ref{robdir}, the Robba ring at $x$ as
$$
\mathcal R_{x} := (i_{x*}\mathcal O_{]x[})_{]\xi[}.
$$
Since $\mathcal O_{X_{K}^{\mathrm{an}},]\xi[}$ is a field, the adjunction map $\mathcal O_{X_{K}^{\mathrm{an}},]\xi[}  \to \mathcal R_{x}$ takes values into the maximum subfield of $\mathcal R_{x}$ which is exactly $\mathcal R^{\mathrm{bd}}_{x}$.

\begin{dfn}
Let $x \in X_{k}$ be a closed point.
If $E$ be a constructible module, the \emph{bounded Robba fiber} of $E$ at  $x$ is
$$
\mathcal R^{\mathrm{bd}}_{x}(E) := \mathcal R_{x}^{\mathrm{bd}} \otimes_{\mathcal O_{X_{K}^{\mathrm{an}},]\xi[}} E_{]\xi[}.
$$
A constructible $\nabla$-module $E$ has \emph{Frobenius type at $x$} if $\mathcal R^{\mathrm{bd}}_{x}(E)$ has a Frobenius structure.
\end{dfn}

Note that this definition is generic in the sense that if $U \subset X_{k}$ is any open subset, then $\mathcal R_{x}^{\mathrm{bd}}(i_{U*}i_{U}^{-1}E) = \mathcal R_{x}^{\mathrm{bd}}(E)$.
We globalize now the definition:

\begin{dfn} \label{bigdef}
A constructible convergent $\nabla$-module $E$ has \emph{Frobenius type} if
\begin{enumerate}
\item $E$ has Frobenius type at all closed points $x \in X_{k}$.
\item $\mathrm R \widetilde{\mathrm{sp}}_{*} E$ is a coherent $\mathcal D^\dagger_{\hat X \mathbf Q}$-module.
\end{enumerate}
\end{dfn}

It is very likely that the second condition is not necessary and it is also possible that the convergence condition is automatic.
Anyway, we have the following:

\begin{prop}
If $F: X \to X$ is a lifting of the absolute Frobenius of $X_{k}$, then any constructible \emph{$F$-$\nabla$-module} on $X_{K}^{\mathrm{an}}$ has Frobenius type.
\end{prop}

\textbf {Proof: }
The Frobenius structure on $E$ will induce a Frobenius structure on all the bounded Robba fibers.
Moreover, we saw in proposition \ref{holon} that $\mathrm R \widetilde{\mathrm{sp}}_{*} E$ is holonomic and in particular coherent.
$\quad \Box$

For computations, it is convenient to have also at our disposal a more algebraic approach.
If $ \mathrm{Spec}\; A \subset X$ is an algebraic lifting of the complement $U$ of $Z$, taking global sections on the map $\mathcal O_{\hat X}({}^\dagger Z)_{\mathbf Q} \to i_{\hat x *} \mathcal O_{\hat x}({}^\dagger x)_{\mathbf Q}$ will give a morphism $A^\dagger_{K} \to \mathcal R^{\mathrm{bd}}_{x}$.
Alternatively, this is the composite map
$$
A^\dagger_{K} \to \mathcal O_{X_{K}^{\mathrm{an}}]\xi[} \to \mathcal R^{\mathrm{bd}}_{x}.
$$

\begin{dfn}
Let $x \in X_{k}$ be a closed point.
Let $U$ be a non empty affine open subset of $X_{k}$ with algebraic lifting $ \mathrm{Spec}\; A$ and $M$ a finite $A^\dagger_{K}$-module. Then, the \emph{bounded Robba fiber} of $M$ at $x$ is
$$
\mathcal R^{\mathrm{bd}}_{x}(M) := \mathcal R_{x}^{\mathrm{bd}} \otimes_{A^\dagger_{K}} M.
$$
A finite $\nabla$-$A^\dagger_{K}$-module $M$ has \emph{Frobenius type at a closed point $x$ (resp. has Frobenius type)}  if the bounded Robba fiber at $x$ (resp. at all closed points $x \in X_{k}$) has a Frobenius structure.

\end{dfn}

Recall that we also introduced in definition \ref{robfib} the Robba fiber of $M$ and that, clearly, $\mathcal R_{x}(M) = \mathcal R_{x} \otimes_{\mathcal R^{\mathrm{bd}}_{x}} \mathcal R^{\mathrm{bd}}_{x}(M)$.

\begin{lem} \label{bdback}
Let $x \in X_{k}$ be a closed point.
\begin{enumerate}
\item Let $Z \subset X_{k}$ be a finite closed subset and $\mathcal E$ be a $\mathcal D^\dagger_{\hat X}({}^\dagger Z)_{\mathbf Q}$-module which is coherent both as $\mathcal D^\dagger_{\hat X \mathbf Q}$-module and as $\mathcal O_{\hat X}({}^\dagger Z)_{\mathbf Q}$-module.
If $M := \Gamma(\hat X, \mathcal E)$, we have
$$
\mathcal R^{\mathrm{bd}}_{x}(\mathcal E) = \mathcal R^{\mathrm{bd}}_{x}(M).
$$
\item Let $U$ be an open subset of $X_{k}$ and $E$ a coherent $i_{U}^{-1}\mathcal O_{X_{K}^{\mathrm{an}}}$-module.
If $M := \Gamma(]U[, E)$, we have
$$
\mathcal R^{\mathrm{bd}}_{x}(E) = \mathcal R^{\mathrm{bd}}_{x}(M).
$$
\end{enumerate}
\end{lem}

\textbf {Proof: }
Since all functors are right exact, the assertions follow from the cases $\mathcal E = \mathcal O_{\hat X}({}^\dagger Z)_{\mathbf Q}$ and $E = i_{U}^{-1}\mathcal O_{X_{K}^{\mathrm{an}}}$.
$\quad \Box$

\begin{prop}
Let $E$ be a constructible convergent $\nabla$-module and $\mathcal E := \mathrm R \widetilde{\mathrm{sp}}_{*} E$.
Assume that $\mathcal E$ is a coherent $\mathcal D^\dagger_{\hat X \mathbf Q}$-module.
Then, if $x \in X_{k}$ is any closed point, we have
$$
\mathcal R^{\mathrm{bd}}_{x}(\mathcal E) = \mathcal R^{\mathrm{bd}}_{x}(E).
$$
In particular $E$ has Frobenius type (at $x$) if and only if $\mathrm R \widetilde{\mathrm{sp}}_{*} E$ has Frobenius type (at $x$).
\end{prop}

\textbf {Proof: }
If $U$ is an open subset of $X_{k}$ with closed complement $Z$, we know that
$$
\mathcal R^{\mathrm{bd}}_{x}(\mathcal E) = \mathcal R^{\mathrm{bd}}_{x}(\mathcal E({}^\dagger Z)) \quad \mathrm{and} \quad \mathrm \mathcal R^{\mathrm{bd}}_{x}(E) = \mathrm \mathcal R^{\mathrm{bd}}_{x}(i_{U*}i_{U}^{-1}E).
$$
We may therefore assume that $E = i_{U*}E_{U}$ with $E_{U}$ coherent as $i_{U}^{-1}\mathcal O_{X^{\mathrm{an}}_{K}}$-module.
And then use lemma \ref{bdback}.
$\quad \Box$

\section{The Deligne-Kashiwara correspondence}

On the analytic side, we will try to stick to the ``connection'' vocabulary and write $\mathrm R \mathcal Hom_{\nabla}(E', E'')$ when $E'$ and $E''$ are two $\nabla$-modules on $X_{K}^{\mathrm{an}}$ for example, but we will systematically identify this space with $\mathrm R \mathcal Hom_{\mathcal D_{X_{K}^{\mathrm{an}}}}(E', E'')$.

\begin{prop} \label{oright}
If $E$ is a $\nabla$-module on $X_{K}^{\mathrm{an}}$, then
$$
\mathrm R\widetilde{\mathrm{sp}}_{*}\mathrm R\mathcal Hom_{\nabla}(\mathcal O_{X_{K}^{\mathrm{an}}}, E) \simeq \mathrm R\mathcal Hom_{\mathcal D_{\hat X\mathbf Q}}(\mathcal O_{\hat X \mathbf Q}, \mathrm R \widetilde{\mathrm{sp}}_{*} E).
$$
Actually, if $E$ is a constructible convergent $\nabla$-module, we have
$$
\mathrm R\widetilde{\mathrm{sp}}_{*}\mathrm R\mathcal Hom_{\nabla}(\mathcal O_{X_{K}^{\mathrm{an}}}, E) \simeq \mathrm R\mathcal Hom_{\mathcal D^\dagger_{\hat X\mathbf Q}}(\mathcal O_{\hat X \mathbf Q}, \mathrm R \widetilde{\mathrm{sp}}_{*} E).
$$
\end{prop}

\textbf {Proof: }
The \emph{analytic Spencer complex}
$$
[\mathcal Hom_{\mathcal O_{X_{K}^{\mathrm{an}}}}(\Omega^1_{X_{K}^{\mathrm{an}}}, \mathcal D_{X_{K}^{\mathrm{an}}}) \to \mathcal D_{X_{K}^{\mathrm{an}}}],
$$
which is locally given by $u \mapsto u(\mathrm d t)\frac \partial {\partial t}$, is a locally free resolution of $\mathcal O_{X_{K}^{\mathrm{an}}}$.
It follows that
$$
\mathrm R\mathcal Hom_{\mathcal D_{X_{K}^{\mathrm{an}}}}(\mathcal O_{X_{K}^{\mathrm{an}}}, E) = E \otimes_{\mathcal O_{X_{K}^{\mathrm{an}}}} \Omega^\bullet_{X_{K}^{\mathrm{an}}}.
$$
But it is also true that the \emph{algebraic Spencer complex}
$$
[\mathcal Hom_{\mathcal O_{\hat X\mathbf Q}}(\Omega^1_{\hat X \mathbf Q}, \mathcal D_{\hat X \mathbf Q}) \to \mathcal D_{\hat X \mathbf Q}]
$$
is a resolution of $\mathcal O_{\hat X \mathbf Q}$ so that
$$
\mathrm R\mathcal Hom_{\mathcal D_{\hat X\mathbf Q}}(\mathcal O_{\hat X \mathbf Q}, \mathrm R \widetilde{\mathrm{sp}}_{*} E) =  [\mathrm R \widetilde{\mathrm{sp}}_{*}E \otimes_{\mathcal O_{\hat X\mathbf Q}} \Omega^\bullet_{\hat X \mathbf Q}].
$$
The first result follows since
$$
\mathrm R \widetilde{\mathrm{sp}}_{*} (E \otimes_{\mathcal O_{X_{K}^{\mathrm{an}}}} \Omega^\bullet_{X_{K}^{\mathrm{an}}}) = [\mathrm R \widetilde{\mathrm{sp}}_{*}E \otimes_{\mathcal O_{\hat X\mathbf Q}} \Omega^\bullet_{\hat X \mathbf Q}].
$$
Now, if $E$ is a constructible convergent $\nabla$-module, we can consider the \emph{arithmetic Spencer complex}
$$
[\mathcal Hom_{\mathcal O_{\hat X\mathbf Q}}(\Omega^1_{\hat X \mathbf Q}, \mathcal D^\dagger_{\hat X \mathbf Q}) \to \mathcal D^\dagger_{\hat X \mathbf Q}],
$$
which is also a resolution of $\mathcal O_{\hat X \mathbf Q}$ (proposition 4.3.3 of \cite{Berthelot00}) and get
$$
\mathrm R\mathcal Hom_{\mathcal D^\dagger_{\hat X\mathbf Q}}(\mathcal O_{\hat X \mathbf Q}, \mathrm R \widetilde{\mathrm{sp}}_{*} E) =  [\mathrm R \widetilde{\mathrm{sp}}_{*}E \otimes_{\mathcal O_{\hat X\mathbf Q}} \Omega^\bullet_{\hat X \mathbf Q}].
\quad \Box$$

\begin{cor} \label{bothco}
If $U$ is an open subset of $X_{k}$ and $E'$ and $E''$ are two coherent $\nabla$-modules on $]U[$, then
$$
\mathrm R\widetilde{\mathrm{sp}}_{*}\mathrm R\mathcal Hom_{\nabla}(i_{U*}E', i_{U*}E'') \simeq \mathrm R\mathcal Hom_{\mathcal D_{\hat X\mathbf Q}}( \widetilde{\mathrm{sp}}_{*} i_{U*}E', \widetilde{\mathrm{sp}}_{*} i_{U*}E'').
$$
Actually, if $E'$ and $E''$ are overconvergent, we have
$$
\mathrm R\widetilde{\mathrm{sp}}_{*}\mathrm R\mathcal Hom_{\nabla}(i_{U*}E', i_{U*}E'') \simeq \mathrm R\mathcal Hom_{\mathcal D^\dagger_{\hat X\mathbf Q}}( \widetilde{\mathrm{sp}}_{*} i_{U*}E', \widetilde{\mathrm{sp}}_{*} i_{U*}E'')
$$
\end{cor}

\textbf {Proof: }
It follows from proposition \ref{spope} that if
$$
E := \mathcal Hom_{i_{U}^{-1}\mathcal O_{X_{K}^{\mathrm{an}}}} (E', E'')
$$
then, we have
$$
E \simeq \mathcal Hom_{\mathcal O_{\hat X}({}^\dagger Z)_{\mathbf Q}}(\widetilde{\mathrm{sp}}_{*}i_{U*}E', \widetilde{\mathrm{sp}}_{*}i_{U*}E'').
$$
One can use the $1$-stratifications in order to show that this isomorphism is horizontal.
Moreover, since $E'$ is coherent, $E$ will be overconvergent when $E'$ and $E''$ are.
And we can then apply the proposition to $i_{U*}E$.
$\quad \Box$

\begin{lem} \label{diif}
Let $U$ be an open subset of $X_{k}$ with closed complement $Z$.
If $E'$ is a locally trivial $\nabla$-module on $]Z[$ and $E''$ is a coherent overconvergent $\nabla$-module on $]U[$, we have
$$
\mathrm R\mathcal Hom_{\nabla}(i_{Z!}E', i_{U*}E'') = 0 \quad \mathrm{and} \quad \mathrm R\mathcal Hom_{\mathcal D^\dagger_{\hat X\mathbf Q}}(\mathrm R \widetilde{\mathrm{sp}}_{*} i_{Z!}E',  \widetilde{\mathrm{sp}}_{*} i_{U*}E'') = 0.
$$
\end{lem}

\textbf {Proof: }
First of all, we have $i_{U}^{-1}i_{Z!}E' = 0$ and therefore, by adjunction,
$$
\mathrm R\mathcal Hom_{\mathcal O_{X_{K}^{\mathrm{an}}}}(i_{Z!}E', i_{U*}E'') = \mathrm R\mathcal Hom_{i_{U}^{-1}\mathcal O_{X_{K}^{\mathrm{an}}}}(i_{U}^{-1}i_{Z!}E', E'') = 0.
$$
The first assertion follows.
For the second one, we can write
$$
\mathrm R \widetilde{\mathrm{sp}}_{*} i_{Z!}E' = i_{Z+}H[-1] \quad \mathrm{and} \quad  \widetilde{\mathrm{sp}}_{*} i_{U*}E'' = \mathcal E
$$
where $H$ is a coherent $\mathcal O_{D_{K}}$-module, $D$ being a smooth lifting of $Z$, and and $\mathcal E$ is a $\mathcal D^\dagger_{\hat X}({}^\dagger Z)_{\mathbf Q}$-module.
Note that
$$
(i_{Z+}H)({}^\dagger Z) := \mathcal D^\dagger_{\hat X}({}^\dagger Z) \otimes_{\mathcal D^\dagger_{\hat X}} i_{Z+}H = 0.
$$
But then, by adjunction, we have
$$
\mathrm R\mathcal Hom_{\mathcal D^\dagger_{\hat X\mathbf Q}}(i_{Z+}H,  \mathcal E) = \mathrm R\mathcal Hom_{\mathcal D^\dagger_{\hat X}({}^\dagger Z)_{\mathbf Q}}((i_{Z+}H)({}^\dagger Z),  \mathcal E) = 0.
\quad \Box$$

Recall that we defined for a constructible module $E$, its Robba fiber (resp. bounded Robba fiber) at a closed point $x \in X_{k}$, as
$$
\mathcal R_{x}(E) = \mathcal R_{x} \otimes_{\mathcal O_{X^{\mathrm{an}}_{K},]\xi[}} E_{]\xi[} \quad \mathrm{(resp.} \quad \mathcal R^{\mathrm{bd}}_{x}(E) = \mathcal R^{\mathrm{bd}}_{x} \otimes_{\mathcal O_{X^{\mathrm{an}}_{K},]\xi[}} E_{]\xi[}).
$$

\begin{lem} \label{dualcom}
Let $U$ be an open subset of $X_{k}$ and $x \not \in U$.
If $E$ is a coherent $i_{U}^{-1}\mathcal O_{X^{\mathrm{an}}_{K}}$-module, then
$$
\mathrm R \mathcal Hom_{\mathcal O_{X_{K}^{\mathrm{an}}}} (i_{U*}E, i_{x!}\mathcal O_{]x[}) = i_{\xi*} \mathrm Hom_{\mathcal R_{x}} (\mathcal R_{x}(E), \mathcal R_{x})
$$
If $E$  is a coherent $\nabla$-module on $]U[$, we have
$$
\mathrm R \mathcal Hom_{\nabla}(i_{U*}E, i_{x!}\mathcal O_{]x[}) = i_{\xi*}  \mathrm Hom_{\nabla} (\mathcal R_{x}(E), \mathcal R_{x}).
$$
\end{lem}

\textbf {Proof: }
By adjunction, we have
$$
\mathrm R \mathcal Hom_{\mathcal O_{X_{K}^{\mathrm{an}}}} (i_{U*}E, i_{x!}\mathcal O_{]x[}) = i_{\xi*} \mathrm Hom_{\mathcal O_{X_{K}^{\mathrm{an}},]\xi[}} (E_{]\xi[}, \mathcal R_{x})
$$
and the first assertion follows by extending scalars.
Then, we can write
$$
\mathrm R \mathcal Hom_{\nabla}(i_{U*}E, i_{x!}\mathcal O_{]x[}) = \mathrm R \mathcal Hom_{\mathcal D_{X_{K}^{\mathrm{an}}}}(\mathcal O_{X_{K}^{\mathrm{an}}}, \mathrm R \mathcal Hom_{\mathcal O_{X_{K}^{\mathrm{an}}}} (i_{U*}E, i_{x!}\mathcal O_{]x[}))
$$
$$
= \mathrm R \mathcal Hom_{\mathcal D_{X_{K}^{\mathrm{an}}}}(\mathcal O_{X_{K}^{\mathrm{an}}}, i_{\xi*} \mathrm Hom_{\mathcal O_{X_{K}^{\mathrm{an}},]\xi[}} (E_{]\xi[}, \mathcal R_{x}))
$$
$$
= i_{\xi*}  \mathrm R Hom_{\mathcal D_{X_{K}^{\mathrm{an}},]\xi[}}(\mathcal O_{X_{K}^{\mathrm{an}},]\xi[},  \mathrm Hom_{\mathcal O_{X_{K}^{\mathrm{an}},]\xi[}} (E_{]\xi[}, \mathcal R_{x})) 
$$
$$
= i_{\xi*} \mathrm {RHom}_{\nabla}(E_{]\xi[}, \mathcal R_{x}) = i_{\xi*}  \mathrm {RHom}_{\nabla}(\mathcal R_{x}(E), \mathcal R_{x})
\quad \Box$$
 
Recall that a constructible $\nabla$-module $E$ is said to have Frobenius type at a closed point $x \in X_{k}$ if its bounded Robba Fiber at this point has a Frobenius structure.

\begin{prop} \label{opp}
Let $E$ be a coherent overconvergent $\nabla$-module on $]U[$ and $x \not \in U$.
Assume $E$ has Frobenius type at $x$ and that $\widetilde{\mathrm{sp}}_{*} i_{U*}E$ is $\mathcal D^\dagger_{\hat X\mathbf Q}$-coherent.
Then, we have
$$
\mathrm R\widetilde{\mathrm{sp}}_{*}\mathrm R\mathcal Hom_{\nabla}(i_{U*}E, i_{x!}\mathcal O_{]x[}) \simeq \mathrm R\mathcal Hom_{\mathcal D^\dagger_{\hat X\mathbf Q}}( \widetilde{\mathrm{sp}}_{*} i_{U*}E,\mathcal H^\dagger_{x})[-1].
$$
\end{prop}

\textbf {Proof: }
Follows from lemma \ref{dualcom} and proposition \ref{crewtech}.
$\quad \Box$

Recall from definition \ref{bigdef}, that a constructible convergent $\nabla$-module $E$ has Frobenius type if it has Frobenius type at all closed points and $\mathrm R \widetilde{\mathrm{sp}}_{*} E$ is a coherent $\mathcal D^\dagger_{\hat X\mathbf Q}$-module.

\begin{thm} \label{rhom}
Let $E'$ and $E''$ be two constructible convergent $\nabla$-module on $X_{K}^{\mathrm{an}}$.
Assume that $E'$ has Frobenius type.
Then, we have
$$
\mathrm R\widetilde{\mathrm{sp}}_{*}\mathrm R\mathcal Hom_{\nabla}(E', E'') \simeq \mathrm R\mathcal Hom_{\mathcal D^\dagger_{\hat X\mathbf Q}}(\mathrm R \widetilde{\mathrm{sp}}_{*} E', \mathrm R \widetilde{\mathrm{sp}}_{*} E'').
$$
\end{thm}

\textbf {Proof: }
Since $E'$ and $E''$ are constructible, there exists an open subset $U$ of $X_{k}$ such that both $i_{U}^{-1}E'$ and $i_{U}^{-1}E''$ are coherent.
If we use corollary \ref{unscrew} for $E''$, it follows from corollary \ref{bothco} and proposition \ref{opp} that
$$
\mathrm R\widetilde{\mathrm{sp}}\;\mathrm R\mathcal Hom_{\nabla}(i_{U*}i_{U}^{-1}E', E'') \simeq \mathrm R\mathcal Hom_{\mathcal D^\dagger_{\hat X\mathbf Q}}( \widetilde{\mathrm{sp}}_{*}i_{U*} i_{U}^{-1}E', \mathrm R \widetilde{\mathrm{sp}}_{*} E'').
$$
If we apply this result to the case $E' = \mathcal O_{X_{K}^{\mathrm{an}}}$ and use corollary \ref{unscrew} for $\mathcal O_{X_{K}^{\mathrm{an}}}$, then proposition \ref{oright} gives
$$
\mathrm R\widetilde{\mathrm{sp}}\;\mathrm R\mathcal Hom_{\nabla}(i_{x!}\mathcal O_{]x[}, E'') \simeq \mathrm R\mathcal Hom_{\mathcal D^\dagger_{\hat X\mathbf Q}}( \mathcal H^\dagger_{x}[-1], \mathrm R \widetilde{\mathrm{sp}}_{*} E'').
$$
And we can apply again corollary \ref{unscrew}, but to $E'$ this time.
$\quad \Box$

\begin{cor} \label{fullyf}
The functor $\mathrm R\widetilde{\mathrm{sp}}_*$ induces a fully faithful functor from the category of constructible convergent $\nabla$-modules of Frobenius type to the category of $\mathcal D^\dagger_{\hat X\mathbf Q}$-modules. $\quad \Box$
\end{cor}

We want now to understand its essential image.
Since we are only interested in curves, we can make the following ad hoc definition:

\begin{dfn}
A coherent $\mathcal D^\dagger_{\hat X\mathbf Q}$-module $\mathcal E$ is \emph{holonomic} if there exists a finite closed subset $Z$ of $X_{k}$ such that $\mathcal E({}^\dagger Z)$ is $\mathcal O_{\hat X}({}^\dagger Z)_{\mathbf Q}$-coherent.
A bounded complex of $\mathcal D^\dagger_{\hat X\mathbf Q}$-modules is \emph{holonomic} if it has holonomic cohomology.
\end{dfn}

Note that an $F$-$\mathcal D^\dagger_{\hat X\mathbf Q}$-module is holonomic in the usual sense if and only if the underlying module is holonomic in this sense.

\begin{prop} \label{essential}
If $\mathcal E$ is a perverse holonomic $\mathcal D^\dagger_{\hat X\mathbf Q}$-module of Frobenius type, there exists a constructible convergent $\nabla$-module $E$ on $X_{K}^{\mathrm{an}}$ such that $\mathrm R\widetilde{\mathrm{sp}}_* E = \mathcal E$.
\end{prop}

\textbf {Proof: }
Let $\mathcal E$ be a perverse holonomic $\mathcal D^\dagger_{\hat X\mathbf Q}$-module.
By definition, it is a complex with holonomic cohomology concentrated in degree 0 and 1, flat in degree 0 and finitely supported in degree 1.
Let $Z$ be a finite closed subset of $X_{k}$ such that $\mathcal H^1(\mathcal E)$ is supported in $Z$ and $\mathcal H^q(\mathcal E)({}^\dagger Z)$ is $\mathcal O_{\hat X}({}^\dagger Z)_{\mathbf Q}$-coherent (for $q = 0, 1$).
Since $\mathcal E$ is coherent, we have an exact triangle
$$
i_{Z+}i_{Z}^{!} \mathcal E \to \mathcal E \to \mathcal E({}^\dagger Z) \to \cdots.
$$

First of all, we see that $\mathcal H^q(\mathcal E({}^\dagger Z)) = \mathcal H^q(\mathcal E)({}^\dagger Z) = 0$ when $q \neq 0$.
In other words, $\mathcal E({}^\dagger Z)$ is a coherent $\mathcal D^\dagger_{\hat X\mathbf Q}$-module which is also $\mathcal O_{\hat X}({}^\dagger Z)_{\mathbf Q}$-coherent.
Thanks to theorem \ref{unpub},  we can write $\mathcal E({}^\dagger Z) = \widetilde{\mathrm{sp}}_* i_{U*}E_{U}$ with $E_{U}$ a coherent overconvergent $\nabla$-module on $]U[$.

By definition, if $x$ is a point in $Z$ and $a$ a non ramified lifting of $x$, we have
$$
i_{x}^{!} \mathcal E = \mathrm Li_{x}^{*} \mathcal E[-1] = 
K(a) \otimes^{\mathrm{L}}_{\mathcal O_{\hat X\mathbf Q,x}} \mathcal E_{x}[-1]
$$
If $x$ is defined as the zero of a local parameter $t$ we can use the free left resolution
$$
\xymatrix{\mathcal O_{\hat X\mathbf Q,x} \ar[r]^t & \mathcal O_{\hat X\mathbf Q,x}}
$$
of $K(a)$ to compute $i_{x}^{!} \mathcal E$.
Then, there is a spectral sequence (or workout the bicomplex if you prefer)
$$
E_{2}^{p,q} := \mathcal H^p\left(\mathcal H^q(\mathcal E_{x}) \stackrel t \to \mathcal H^q(\mathcal E_{x})\right)\Rightarrow i_{x}^{!} \mathcal E
$$
with $H^q(\mathcal E_{x})$ flat in degree 0, $t$-torsion in degree 1 and 0 otherwise.
It follows that $i_{x}^{!} \mathcal E$ is just a finite dimensional vector space placed in degree 1.
Thanks to proposition \ref{fineq}, we can therefore write $i_{Z+}i_{Z}^{!} \mathcal E = \mathrm R\widetilde{\mathrm{sp}}_* i_{Z!}E_{Z}$ with $E_{Z}$ a locally trivial $\nabla$-module on $]Z[$.
And the assertion then follows from theorem \ref{rhom}.
$\quad \Box$

\begin{cor}
The functor $\mathrm R\widetilde{\mathrm{sp}}_*$ induces an equivalence between the category of constructible convergent $\nabla$-modules of Frobenius type on $X_{K}^{\mathrm{an}}$ and perverse holonomic $\mathcal D^\dagger_{\hat X\mathbf Q}$-modules of Frobenius type.
$\quad \Box$
\end{cor}

We can at last prove the overconvergent Deligne-Kashiwara correspondence for curves (when there exists a global lifting of Frobenius):

\begin{thm}
Assume that $F: X \to X$ is a lifting of the absolute Frobenius of $X_{k}$.
Then, the functor $\mathrm R\widetilde{\mathrm{sp}}_*$ induces an equivalence between constructible $F$-$\nabla$-modules on $X_{K}^{\mathrm{an}}$ and perverse holonomic $F$-$\mathcal D^\dagger_{\hat X\mathbf Q}$-modules on $\hat X$.
\end{thm}

\textbf {Proof: }
Follows from corollary \ref{fullyf}, proposition \ref{holon} and proposition \ref{essential}. 
$\quad \Box$

Note that this theorem extends to the case where there is no global lifting of Frobenius.
Actually, one expects an analog result in higher dimension.
However, we believe that our methods reach their limit here and that a crystalline approach will be more appropriate to go further.

\bibliographystyle{plain}
\addcontentsline{toc}{section}{Bibliography}
\bibliography{BiblioBLS}

\end{document}